\documentclass[11pt]{amsart}
\usepackage[utf8]{inputenc}
\usepackage{amsmath,amsthm,amsfonts,amssymb,amscd,mathtools}
\usepackage{verbatim}
\usepackage{mathrsfs}
\usepackage{multicol}
\usepackage{tabularx}
\usepackage[usenames,dvipsnames]{color}
\usepackage{enumerate}
\usepackage{graphicx}
\usepackage{color,xcolor}
\usepackage{bbm}
\usepackage[all]{xy}
\usepackage{hyperref}
\usepackage[capitalise]{cleveref}
\usepackage{pb-diagram,tikz-cd}

\theoremstyle{definition}
\newtheorem{thm}{Theorem}[subsection]
\newtheorem{prop}[thm]{Proposition}
\newtheorem{cor}[thm]{Corollary}

\newtheorem{lem}[thm]{Lemma}
\newtheorem{defn}[thm]{Definition}
\newtheorem{ex}[thm]{Example}
\newtheorem{rem}[thm]{Remark}

\numberwithin{equation}{subsection}

\def\bs#1{\boldsymbol{#1}}

\def\lie#1{\mathfrak{#1}}
\def\tlie#1{\tilde{\mathfrak{#1}}}

\def\endd{\hfill$\diamond$}

\DeclareMathOperator{\Adj}{Adj}

\begin{document}

	\title[Reality Determining Subgraphs and Strongly Real Modules]{Reality Determining Subgraphs\\ and Strongly Real Modules}
	
	\author{Matheus Brito}
	\address{Departamento de Matemática, Universidade Federal do Paraná, Curitiba - PR - Brazil, 81530-015}
	\email{mbrito@ufpr.br, cristiano.clayton@ufpr.br}
	\thanks{}
	\author{Adriano Moura}
	\address{Departamento de Matemática, Universidade Estadual de Campinas, Campinas - SP - Brazil, 13083-859.}	\email{aamoura@unicamp.br}
	\thanks{The work of M.B. and A.M. was partially supported CNPq grant 405793/2023-5. The work of A.M. was also partially supported by Fapesp grant 2018/23690-6.  Part of the ideas of this paper were developed while he visited the Schenzhen International Center for Mathematics. He would like to thank V. Futorny for the opportunity to visit the Center, as well as all the staff for the support leading to excellent working conditions.}
	\author{Clayton Silva}
	\thanks{C. Silva is grateful to CNPq for financial support (grant 100161/2024-3 INCTMat-IMPA) during his postdoctoral research internship at UFPR}
	
	\begin{abstract}
		The concept of pseudo $q$-factorization graphs was recently introduced by the last two authors as a combinatorial language which is suited for capturing certain properties of Drinfeld polynomials. Using certain known representation theoretic facts about tensor products of Kirillov Reshetikhin modules and $q$-characters, combined with special topological/combinatorial properties of the underlying $q$-factorization graphs, the last two authors showed that, for algebras of type $A$, modules associated to totally ordered graphs are prime, while those associated to trees are real. In this paper, we extend the latter result. We introduce the notions of strongly real modules and that of trees of modules satisfying certain properties. In particular, we can consider snake trees, i.e., trees formed from snake modules. Among other results, we show that a certain class of such generalized trees, which properly contains the snake trees, give rise to strongly real modules.
	\end{abstract}
	
	\maketitle

\section{Introduction}

The monoidal structure of the category $\mathcal C$ of finite-dimensional modules for a quantum affine algebra $U_q(\tlie g)$ has been an intensive topic of research from a variety of perspectives. The original motivation for the study was the connection with branches of Mathematical Physics. However, the rich underlying combinatorics lead to interactions with the theory of cluster algebras \cite{HL:cluster}, Grassmanians \cite{cdf:grass}, representations of $p$-adic groups \cite{LM:sqirred}, tropical geometry \cite{eali:trop}, etc. 

Possibly the most natural problem when working with an abelian monoidal category, such as $\mathcal C$, is the description of the structure of the tensor product of two given simple objects. There is a class of simple modules in $\mathcal C$ which can be regarded as building blocks of the monoidal structure: that of Kirillov-Reshetikhin (KR) modules. Recall the simple modules in $\mathcal C$ are highest-$\ell$-weight (i.e., they are highest-weight with respect to a certain ``triangular decomposition'') and the corresponding highest $\ell$-weight can be encoded in an $I$-tuple of polynomials with constant term $1$ known as Drinfeld polynomials, where $I$ is the set of nodes of the Dynkin diagram of the underlying simple Lie algebra $\lie g$. We denote by $V(\bs\pi)$ a simple module with Drinfeld polynomial $\bs\pi$. The Drinfeld polynomial of a KR module has only one non-constant entry and the roots of the non-constant entry form a $q^{d_i}$-string (the positive integers $d_i$ symmetrize the Cartan matrix of $\lie g$). We denote such tuples by $\bs\omega_{i,a,r}$, where $i$ is the node corresponding to the non-constant entry, $a$ and $r$ are the center and the length of the $q$-string, respectively. It is known that there exist finite sets of positive integers $\mathscr R_{i,j}^{r,s}$ such that
\begin{equation*}
	V(\bs\omega_{i,a,r})\otimes V(\bs\omega_{j,b,s}) \ \ \text{is reducible} \quad\Leftrightarrow\quad \frac{a}{b} = q^m \ \ \text{with}\ \ |m|\in\mathscr R_{i,j}^{r,s}.
\end{equation*}
Moreover, when such tensor product is reducible, it is highest-$\ell$-weight if and only if $m> 0$.   Beyond the class of KR modules, there is no explicit general answer for determining if a tenor product of two simple modules is simple or not. We let $\mathcal{KR}$ denote the subset of KR type Drinfeld polynomials inside the set $\mathcal P^+$ of all  Drinfeld polynomials. 

This paper focuses on the question of determining whether the tensor square of a given simple module is simple. Such modules are called real due to the relation with cluster algebras.	In particular, since $0\notin\mathscr R_{i,i}^{r,r}$, KR modules are real.  In \cite{hl:KR}, Hernandez and Leclerc introduced a proper subcategory $\mathcal C^-$ of $\mathcal C$ which ``essentially'' contains all simple objects. They proved an
algorithm for calculating $q$-characters of KR modules by studying the cluster algebra $\mathscr A$ which is
isomorphic to the Grothendieck ring of the subcategory $\mathcal C^-$. The authors also conjectured that all cluster monomials of $\mathscr A$ correspond to the classes of real 
simple objects of $\mathcal C^-$ (see \cite[Conjecture 5.2]{hl:KR}). 	Kashiwara, Kim, Oh, and Park proved this conjecture by studying a complete duality datum which induces a Schur--Weyl duality functor from the category of finite dimensional graded modules for some KLR-algebra to $\mathcal C^-$ \cite{kkop:mc1, KKOP:cluster, kkop:mc2}. Prior to that, in \cite{Qin}, Qin proved this conjecture for an infinity family of subcategories of $\mathcal C^-$ when $\lie g$ is simply laced, using different methods which rely on the geometric approach of quiver varieties. 

 Thus, in principle, the Drinfeld polynomials of simple real modules can be determined using cluster algebra machinery. In practice, determining if a given Drinfeld polynomial is a cluster monomial is not a simple task. One particularly important class of modules for which the connection with cluster monomials can be made explicit is that of snake modules introduced in \cite{my:pathB}, which include the KR modules and, more generally, the minimal affinizations. In \cite{DLL:snakes}, the authors introduce $S$-systems consisting of equations satisfied by the $q$-characters of prime snake modules of types $A$ and $B$. They also show that every equation in the $S$-system corresponds to a mutation in the associated cluster algebra and every prime snake module corresponds to some cluster variable.

In this paper, we approach the task from a different perspective: can we create a combinatorial apparatus which would lead to a ``simple'' description of classes of Drinfeld polynomials of real modules? For that, we use the notion of pseudo $q$-factorization graph introduced in \cite{ms:to}. Quite clearly, any Drinfeld polynomial can be written as a product of KR type Drinfeld polynomials. We call any such factorization a pseudo $q$-factorization. The actual $q$-factorization is the one with factors of ``maximal degree'', while the fundamental factorization is that whose all factors have degree $1$. The multiset of (pseudo) $q$-factors is then regarded as the vertex set $\mathcal V$ of a digraph and a pair of vertices $(\bs\omega,\bs\varpi)$ is an arrow if and only if
\begin{equation*}
	V(\bs\omega)\otimes V(\bs\varpi) \quad\text{is reducible and highest-$\ell$-weight.}
\end{equation*} 
The graph $G$ is then called a (pseudo) $q$-factorization graph over the given $\bs\pi\in\mathcal P^+$. Any subgraph $H$ of $G$ gives rise to an element $\bs\pi_H\in\mathcal P^+$ which divides $\bs\pi$ in the obvious manner. 

Given two $q$-factorization graphs $G_1$ and $G_2$, we let $G_1\otimes G_2$ denote the graph whose vertex set is the union of those of $G_1$ and $G_2$. In particular, the pair $(G_1,G_2)$ is a cut of the graph $G_1\otimes G_2$. Suppose $G_1$ and $G_2$ are real, i.e., $V(\bs\pi_{G_k})$ is a real module. Under which conditions $G_1\otimes G_2$ is also real? A sufficient condition for this is given by \Cref{t:isreal}. This theorem has similar flavor and shares some assumptions with that of \cite[Lemma 2.10]{LM:sqirred}, although they are not exactly the same. Equivalently, given a graph $G$, can we find a subgraph $H$ such that the corresponding cut $(H,H^c)$, where $H^c$ denote the subgraph complementary to $H$ in $G$, satisfy the requirements of \Cref{t:isreal}? This leads to one of the main definitions of the paper: that of reality determining subgraphs (rds for short) -- see \Cref{d:rds}. In particular, $G$ is an rds for itself if, and only if, it has a single vertex. \Cref{d:rds} is designed so that it follows from \Cref{t:isreal} that $G$ is real if it contains an rds.

The main constructions of this paper arise from combining the notion of rds with a purely graph theoretical concept which we termed ``quochains''. The choice of terminology arises from the interpretation of $H^c$ as a quotient of $G$ by the subgraph $H$. By a multicut of $G$,  we mean a sequence $\mathcal G=G_1,\dots,G_l$ of subraphs with pairwise disjoint vertex sets whose union is the vertex set of $G$. Thus, $G_2$ can be regarded as a subgraph of the quotient $G/G_1$ of $G$ by $G_1$, $G_3$ as a subgraph of the quotient  $(G/G_1)/G_2$, and so on. Thus, we have a chain of subgraphs of quotients -  a quochain. In terms of the tensor product notation mentioned above, if we set $\bar G_{k-1} = G_k\otimes\cdots\otimes G_l$ for $1\le k\le l$, then $G_k$ is interpreted as a subgraph of $\bar G_{k-1}$. By an rds-quochain, we mean a quochain such that $G_k$ is an rds in $\bar G_{k-1}$ for all $1\le k\le l$. In particular, $G_l$ must be a singleton.   

Finally, the notion of rds-quochains lead to that of strongly real graphs and modules. We say a graph $G$ is strongly real if it admits an rds-quochain such that all parts of the underlying multicut are singletons, i.e., correspond to KR modules. A simple module is said to be strongly real if it affords a strongly real pseudo $q$-factorization graph. In other words, if its Drinfeld polynomial can be built by a sequence of KR type Drinfeld polynomials such that every partial product of this sequence corresponds to a real module. In Examples \ref{ex:realbutnotstrong} and \ref{ex:realbutnotstrongkkop}, we present real modules which are not strongly real, but whose reality can be detected by an rds. We found no example of a real graph which admits no proper rds (see the discussion in \Cref{ss:nords}, where we explore the notion of weak rds-quochain).   

In Sections \ref{ss:eas} and \ref{ss:etoas}, we describe sufficient conditions for a vertex to be an rds. As a corollary, we recover one of the main results of \cite{ms:3tree} which says that trees are real (in fact, strongly real) and that prime snake modules are real (in fact, strongly real). In \Cref{ss:srt}, we further explore the notion of multicuts to introduce a generalization of the notion of tree in graph theory. Namely, given a sequence $\mathcal G=G_1,\dots, G_l$ of graphs, we say the graph $G=G_1\otimes\cdots\otimes G_l$ is a $\mathcal G$-tree if the number of arrows in the corresponding cut set is $l-1$. The cut set is the set of arrows of $G$ which are not arrows of any of the parts $G_k$ of the multicut. We then address the problem of determining if a graph which can be realized as a $\mathcal G$-tree whose parts admit (weak) rds-quochains, also admit (weak) rds-quochains. The main results related to this are Theorems \ref{t:srtree} and \ref{t:kkoprds}. For the proofs, we systematically use certain invariants introduced in \cite{kkop:mc1}, which we refer to as the KKOP invariants. 

It is interesting to emphasize that, although Theorems \ref{t:srtree} and \ref{t:kkoprds} explore the same notions, they have somewhat opposite philosophies for forming the $\mathcal G$-trees. In \Cref{t:srtree}, we request the parts of the multicut $\mathcal G$ have very strong properties related to the concept of rds, but we impose no restriction on the KKOP invariants of the pairs of parts. For instance, in \Cref{ss:etoas}, we prove that any vertex of a pseudo $q$-factorization graph over a prime snake module is the final vertex of a strong rds-quochain. In particular, \Cref{t:srtree} applies to the class of snake trees, i.e., $\mathcal G$-trees whose parts are graphs over prime snake modules, and it follows that snake trees are strongly real. We give examples of such modules which are not snake modules, nor their graphs are trees. It would be interesting to study how these modules can be interpreted from the point of view of cluster algebras, Grassmanians, tropical geometry, p-adic groups, etc. We leave this task for forthcoming publications. On the other hand, in \Cref{t:kkoprds}, we impose no restriction on the parts of the multicut except that the underlying modules are real, but request the associated KKOP invariants are at most $1$. In that case, the theorem says that $\mathcal G$ is itself a weak rds quochain, up to reordering.

In \Cref{ex:srnstkkop}, we study a module which does not fit in the assumptions of Theorems \ref{t:srtree} and \ref{t:kkoprds}, but a suitable combination of these theorems implies it is strongly real. The examples of strongly real modules we give here belong to the class of modules realizable as $\mathcal G$-trees for certain proper multicuts of a $q$-factorization $G$, except for the prime snake modules themselves. The case of prime snake modules show the methods can be used in a broader context, in principle.

In \cite{LM:sqirred}, the real simple objects of $\mathcal C$ satisfying a certain condition called regular have been classified when $\lie g$ is of type $A$  (in the 	language of representations of $p$-adic groups). The notion of regularity is related to the classification of rationally smooth Schubert varieties in type $A$ flag varieties. The aforementioned \Cref{ex:realbutnotstrong} is an example of a real module which does not satisfy the regularity condition. Although it is not strongly real, it is ``the next best thing'' (from the perspective of the discussion of \Cref{ss:nords}). Although the discussion about this regularity condition escapes our scope of expertise, it would be interesting to understand examples such as \ref{ex:realbutnotstrong} and \ref{ex:realbutnotstrongkkop} from that perspective. For instance: Are these examples regularizable via Zelevinsky involution? Or is the fact that they are not strongly real an obstruction for regularization? 

Let us comment about how our results depend on the Cartan type of $\lie g$. None of the arguments developed here depend on the type. However, previously proved results with type restrictions are used at some crucial moments. We believe all of them have at least partial extensions for general $\lie g$ and we make local remarks about this. 

The paper is organized as follows. The preliminary background is reviewed in \Cref{s:baseback}, while the notion of pseudo $q$-factorization graphs is reviewed in \Cref{ss:qfgraphs}. The main definitions of the paper are given in \Cref{ss:realtest}, motivated by the aforementioned \Cref{t:isreal}. In the analysis of \Cref{ex:realbutnotstrong}, and also of \Cref{ex:srnst} later on,  we use one of the main results of \cite{ms:3tree} about prime factorization of modules with $3$ $q$-factors for type $A$, which is reviewed in \Cref{ss:3v}. 
The contents of Sections \ref{ss:eas} and \ref{ss:etoas} have already been described above. \Cref{s:applications} is dedicated to applying the newly introduced concepts to describe, in the language of graphs, interesting examples of real modules. The notion of multicuts of tree type is introduced in \Cref{ss:srt}, while Sections \ref{ss:snaketrees} and \ref{ss:kkkords} are dedicated to describing different ways of exploring all of the concepts to generate examples of (strongly) real modules, synthetized by the aforementioned Theorems \ref{t:srtree} and \ref{t:kkoprds}. 
We end the paper with two sections with discussions about topics for future analysis. In \Cref{ss:cluster}, we make some comments about the relation with cluster algebras, while in \Cref{ss:nords} we propose a stratification of the set of (isomorphism classes of) real modules by certain indexes related to \Cref{d:weakrdsq}. \\

\noindent{\em Acknowledgments. The authors thank the referee for their
careful reading of the paper and the many helpful comments and remarks.}

\section{Basic Notation and Background}\label{s:baseback}
	
Throughout the paper, let  $\mathbb Z$ denote the  set integers. Let also $\mathbb Z_{\ge m} ,\mathbb Z_{< m}$, etc.,  denote the obvious subsets of $\mathbb Z$. Given a ring $\mathbb A$, the underlying multiplicative group of units is denoted by $\mathbb A^\times$. 
The symbol $\cong$ means ``isomorphic to''. We shall use the symbol $\diamond$ to mark the end of remarks, examples, and statements of results whose proofs are postponed. The symbol \qedsymbol\ will mark the end of proofs as well as of statements whose proofs are omitted. 
	
\subsection{Cartan Data and Quantum Algebras}\label{ss:basenot}
Let $\lie g$ be a finite-dimensional simple Lie algebra over $\mathbb C$ of rank $n$ and let $I$ be the set of nodes of its Dynkin diagram. We let $x_i^\pm, h_i, i\in I$, denote generators as in Serre's Theorem and let $\lie g=\lie n^-\oplus\lie h\oplus \lie n^+$ be the corresponding triangular decomposition. 
The weight lattice of $\lie g$ will be denoted by $P$ and $P^+$ denotes the subset of dominant weights. The root system, set of positive roots, root lattice and the corresponding positive cone will be denoted, respectively, by $R, R^+, Q,Q^+$, while the fundamental weights and simple roots will be denoted by $\omega_i, \alpha_i,i\in I$.   For $i\in I$, let $i^*= w_0(i)$, where $w_0$ is the Dynkin diagram automorphism induced by the longest element of the Weyl group.  Fix also relatively prime positive integers $d_i, i\in I$, such that $DC$ is symmetric, where $D={\rm diag}(d_1,\dots,d_n)$ and $C$ is the Cartan matrix of $\lie g$, i.e., $c_{i,j}=\alpha_j(h_i)$. 

Consider the quantum affine (actually loop) algebra $U_q(\tlie g)$ over an algebraically closed field of characteristic zero $\mathbb F$, where $q\in\mathbb F^\times$ is not a root of unity.   We use the presentation in terms of generators and relations and the Hopf algebra structure as in  \cite{Moura}. In particular, the generators are denoted by $x_{i,r}^\pm, k_i^{\pm 1}, h_{i,s}, i\in I, r,s\in\mathbb Z, s\ne 0$. The subalgebra generated by $x_i:=x_{i,0}^\pm, k_i^{\pm 1}, i\in I$ is a Hopf subalgebra of $U_q(\tlie g)$ isomorphic to the quantum algebra $U_q(\lie g)$. 
We let $\bs\kappa$ denote the Cartan involution, i.e., the unique algebra automorphism of $U_q(\tlie g)$ such that (see \cite[Proposition 3.2.1(c)]{ms:to} and references therein)
\begin{equation*}
	\kappa(x_{i,r}^{\pm})=-x_{i,-r}^{\mp},\quad \kappa(h_{i,s})=-h_{i,-s},\quad \kappa(k_{i}^{\pm1})=k_{i}^{\mp1}, \quad i\in I,\ r,s\in\mathbb Z, s\ne 0.
\end{equation*}

If $J\subseteq I$, we let $U_q(\tlie g)_J$ denote the subalgebra generated by these elements with $i\in J$. It is naturally isomorphic to the algebra  $U_{q_J}(\tlie g_J)$ associated to the diagram subalgebra $\lie g_J$ of $\lie g$, where $q_J=q^{d_J}$ with $d_J = \min\{d_j:j\in J\}$, but it is not a Hopf subalgebra. In particular, if $J=\{i\}$ for some $i\in I$, then $U_q(\tlie g)_J\cong U_{q_i}(\tlie{sl}_2)$, where $q_i=q^{d_i}$.  We denote by $\check h$ the dual Coxeter number of $\lie g$ and by $\check h_J$ the one of $\lie g_J$.
	
\subsection{Finite-Dimensional Modules}\label{ss:fdreps}
For  $i\in I$, $a\in \mathbb Z$, we let $\bs\omega_{i,a}$ denote the corresponding fundamental $\ell$-weight, which is the Drinfeld polynomial whose unique non-constant entry is equal to $1-q_i^au\in\mathbb F[u]$. 
We let $\mathcal P^+$ denote the monoid multiplicatively generated by such elements and let $\mathcal P$ be the corresponding abelian group. The identity element will be denoted by $\bs 1$. We shall also use an alternative notation for the fundamental $\ell$-weights which is often better suited for the visualization  of the graphs below. Namely, we write $i_a$ for $\bs\omega_{i,a}$.

Let $\mathcal C$ be the full subcategory of that of fintie-dimensional $U_q(\tlie g)$-modules whose simple factors have highest $\ell$-weights in $\mathcal P^+$ and, hence, $\ell$-weights in $\mathcal P$. 
For $\bs\pi\in\mathcal P^+$, $V(\bs\pi)$ will denote a simple $U_q(\tlie g)$-module whose highest $\ell$-weight is $\bs\pi$.
In particular, $V(\bs 1)$ is the trivial one-dimensional module.
$\mathcal C$ is a monoidal abelian category and, hence, the notion of prime objects is defined. Since $V(\bs 1)$ is the unique invertible object in $\mathcal C$, then $V\in\mathcal C$ is prime if 
\begin{equation*}
	V\cong V_1\otimes V_2 \quad\Rightarrow\quad V_j\cong V(\bs 1) \quad\text{for some}\quad j\in\{1,2\}.
\end{equation*}
This definition says $V(\bs 1)$ is prime. Every simple module admits a decomposition as a tensor product of simple prime modules. A simple object $V$ is said to be real if $V\otimes V$ is simple.

For an object $V\in\mathcal C$, let $V^*$ and ${}^*V$ be the dual modules to $V$ such that the usual evaluation maps
\begin{equation*}
	V^*\otimes V \to \mathbb F  \quad\text{and}\quad V\otimes {}^*V\to\mathbb F
\end{equation*}
are module homomorphisms (cf. \cite[Section 2.6]{ms:to}). Then $({}^*V)^*\cong V\cong {}^*(V^*)$ and 
\begin{equation}\label{e:dualtp}
	(V_1\otimes V_2)^*\cong V_2^*\otimes V_1^*.
\end{equation}
The Hopf algebra structure on $U_q(\tlie g)$ is chosen so that, if $V=V(\bs\pi)$, then $V^*\cong V(\bs\pi^*)$, where $\bs\pi\mapsto \bs\pi^*$ is the group automorphism of $\mathcal P$ determined by
\begin{equation*}
	\bs\omega_{i,a}^* = \bs\omega_{i^*,a-\check h}. 
\end{equation*}
Similarly, considering the automorphism determined by ${}^*\bs\omega_{i,a} = \bs\omega_{i^*,a+\check h}$, it follows that ${}^*V(\bs\pi)\cong V({}^*\bs\pi)$. 
Denote by $V^\kappa$ the module obtained from $V$ by pulling back the action through the Cartan involution. Then, $V(\bs\pi)^\kappa\cong V(\bs\pi^\kappa)$ where $\bs\pi\mapsto \bs\pi^\kappa$ is the group automorphism determined by 
\begin{equation*}
	\bs\omega_{i,a}^\kappa = (\bs\omega_{i,-a})^*.
\end{equation*}
In particular,
\begin{equation}\label{e:cart-}
	{}^*(\bs\omega_{i,a}^\kappa) = \bs\omega_{i,-a}.
\end{equation}
	
If $J\subseteq I$, we shall denote by $V(\bs\pi)_J$ the $U_q(\tlie g)_J$-submodule of $V(\bs\pi)$ generated by the top weight space. Under the natural isomorphism $U_q(\tlie g)_J\cong U_q(\tlie g_J)$, $V(\bs\pi)_J$ is isomorphic to $V(\bs\pi_J)$, where $\bs\pi_J$ is the associated $J$-tuple of polynomials. We shall abuse notation and regard $V(\bs\pi)_J$ both as a $U_q(\tlie g)_J$-module and as a $U_q(\tlie g_J)$-module as, typically, no contextual confusion will arise. More generally, if $V$ is a highest-$\ell$-weight module with highest-$\ell$-weight vector $v$, we let $V_J$ denote the $U_q(\tlie g)_J$-submodule of $V(\bs{\pi})$ generated by $v$. Evidently, if $\bs\pi$ is the highest-$\ell$-weight of $V$, then $V_J$ is highest-$\ell$-weight with highest $\ell$-weight $\bs\pi_J$.	
	
Given $i\in I, a\in\mathbb Z, m\in\mathbb Z_{\ge 0}$,  define 
\begin{equation*}
	\bs\omega_{i,a,r} = \prod_{p=0}^{r-1} \bs\omega_{i,a+{r-1-2p}}.
\end{equation*}
We shall refer to Drinfeld polynomials of the form $\bs\omega_{i,a,r}$ as polynomials of Kirillov-Reshetikhin (KR) type, since the corresponding simple modules $V(\bs\omega_{i,a,r})$ are known as Kirillov-Reshetikhin  modules.  The number $a$ will be referred to as the center of the associated $q_i$-string.  The set of all such polynomials will be denoted by $\mathcal{KR}$. 
Every $\bs\pi\in\mathcal P^+$ can be written uniquely as a product of KR type polynomials such that, for every two factors supported at $i$, say $\bs\omega_{i,a,r}$ and $\bs\omega_{i,b,s}$, the following holds
\begin{equation}\label{e:defqfact}
	a-b \notin\mathscr R^{r,s} :=\{r+s-2p: 0\le p<\min\{r,s\}\}. 
\end{equation}
Such factorization is said to be the  $q$-factorization of $\bs\pi$ and the corresponding factors are called the $q$-factors of $\bs\pi$.
By abuse of language, whenever we mention the set of $q$-factors of $\bs\pi$ we actually mean the associated multiset of $q$-factors counted with multiplicities in the $q$-factorization.  It is often convenient to work with factorizations in KR type polynomials which not necessarily satisfy \eqref{e:defqfact}. Such a factorization will be referred to as a pseudo $q$-factorization and the associated factors as the corresponding $q$-factors of the factorization. The set in \eqref{e:defqfact} determines the reducibility of tensor products of KR modules for type $A_1$-subdiagrams as explained in \eqref{e:redsl2sd} below.

\subsection{Tensor Products}\label{ss:tpcrit} We now collect a few known facts about tensor products of highest-$\ell$-weight modules. 
	The following is well-known.
	
    \begin{prop}\label{p:sinter}
		Let $\bs{\pi}, \bs{\varpi}\in\mathcal{P}^+$. Then, $V(\bs{\pi})\otimes V(\bs{\varpi})$ is simple if and only if $V(\bs{\varpi})\otimes V(\bs{\pi})$ is simple and, in that case, $V(\bs{\pi})\otimes V(\bs{\varpi}) \cong V(\bs\pi\bs\varpi)\cong V(\bs{\varpi})\otimes V(\bs{\pi})$. 	\hfil\qed
	\end{prop} 
	
	We also have:
	
	\begin{prop}[{\cite[Corollary 4.1.4]{ms:to}}]\label{p:vnvstar}
		Let $\bs{\pi}, \bs\varpi\in\mathcal{P}^+$. Then, $V(\bs{\pi})\otimes V(\bs\varpi)$ is simple if and only if both $V(\bs{\pi})\otimes V(\bs\varpi)$ and $V(\bs\varpi)\otimes V(\bs{\pi})$ are highest-$\ell$-weight. 	\hfil\qed
	\end{prop}
	
Given a connected subdiagram $J$, since $U_q(\tlie g)_J$ is not a sub-coalgebra of $U_q(\tlie g)$, if $M$ and $N$ are $U_q(\tlie g)_J$-submodules of $U_q(\tlie g)$-modules $V$ and $W$, respectively, it is in general not true that $M\otimes N$ is a $U_q(\tlie g)_J$-submodule of $V\otimes W$.  Recalling that we have an algebra isomorphism $U_q(\tlie g)_J\cong U_{q_J}(\tlie g_J)$, denote by $M\otimes_J N$ the $U_q(\tlie g)_J$-module obtained by using the coalgebra structure from $ U_{q_J}(\tlie g_J)$. The next result describes a special situation on which $M\otimes N$ is a submodule isomorphic to $M\otimes_J N$. It is known that, if $V$ and $W$ are finite-dimensional highest-$\ell$-weight modules, then $V_J\otimes W_J$ is a $U_q(\tlie g)_J$-submodule of $V\otimes W$ and the identity map induces an isomorphism of  $U_q(\tlie g)_J$-modules
\begin{equation}\label{e:subJtp}
		V_J\otimes W_J \cong V_J\otimes_J W_J.
\end{equation}
Moreover (see \cite[Corollary 3.2.4]{ms:to} and references therein),
\begin{equation}\label{e:sJs}
	V\otimes W \ \ \text{highest-$\ell$-weight (simple)} \quad\Rightarrow\quad  V_J\otimes W_J\ \ \text{highest-$\ell$-weight (simple)}.
\end{equation}

The next two proposition were among the main tools used in \cite{ms:to,ms:3tree} to prove certain tensor products are highest-$\ell$-weight.
	
\begin{prop}[{\cite[Propositions 4.3.1]{ms:to}}]\label{p:hlwmorph}
	Let $\bs{\omega},\bs{\varpi}\in\mathcal{P}^+$ and  $V=V(\bs{\omega})\otimes V(\bs{\varpi})$. Then, $V$ is highest-$\ell$-weight provided there exists  $\bs{\mu}\in\mathcal{P}^+$ such that one of the following conditions holds:
\begin{enumerate}[(i)]
	\item $V(\bs{\omega}\bs{\mu})\otimes V(\bs{\varpi})$ and $V(\bs{\omega})\otimes V(\bs{\mu})$ are both highest-$\ell$-weight;
	\item $V(\bs{\omega})\otimes V(\bs{\mu}\bs{\varpi})$ and $V(\bs{\mu})\otimes V(\bs{\varpi})$ are both highest-$\ell$-weight. 	\qed
\end{enumerate}
\end{prop}

\begin{prop}[{\cite[Propositions 4.5.1]{ms:to}}]\label{p:hlwdualmorph}
	Let $\bs{\lambda},\bs{\mu},\bs{\nu}\in\mathcal{P}^+$. Let also $V = V(\bs{\lambda})\otimes V(\bs{\nu})^*$,
	\begin{equation*}
		T_1=V(\bs{\lambda}\bs{\mu})\otimes V(\bs{\nu}),\quad U_1=V(\bs{\lambda})\otimes V(\bs{\mu}), \quad W_1=V(\bs{\mu})\otimes V(\bs{\nu}),
	\end{equation*}	
	\begin{equation*}
		T_2=V(\bs{\lambda})\otimes V(\bs{\mu}\bs{\nu}),\quad U_2=V(\bs{\mu})\otimes V(\bs{\nu}),\quad W_2=V(\bs{\lambda})\otimes V(\bs{\mu}).
	\end{equation*}
	Then, $W_i$ is highest-$\ell$-weight provided $T_i$ and $U_i$ are highest-$\ell$-weight,  $i\in\{1,2\}$, and $V$ is simple.\qed
\end{prop}

The following theorem plays a crucial role in the main results of this paper. For comments on its proof, see \cite[Remark 4.1.7]{ms:to}.
The statement we reproduce here is essentially that of \cite[Theorem 4.2]{gumi}.
	
\begin{thm}\label{t:cyc}
	Let $S_1,\cdots, S_m\in\mathcal C$ be simple and assume $S_i$ is real either for all $i>2$ or for all $i<m-1$. If $S_i\otimes S_j$ is highest-$\ell$-weight for all $1\leq i< j\leq m$, then $S_1\otimes\cdots\otimes S_m$ is highest-$\ell$-weight. 	\hfil\qed
\end{thm}
	
\begin{cor}\label{c:hlwquot}
	Let $\bs{\pi},\bs\varpi\in\mathcal{P}^+$ and suppose  there exist 
	$\bs\pi_k,\bs\varpi_l\in\mathcal P^+, 1\le k\le m, 1\le l\le m'$, such that  
	\begin{equation*}
		\bs\pi = \prod_{k=1}^m \bs\pi_k, \quad  \bs\varpi = \prod_{k=1}^{m'} \bs\varpi_k, 
	\end{equation*}
	and the following tensor products are highest-$\ell$-weight:
	\begin{equation*}
		V(\bs\pi_k)\otimes V(\bs\pi_l), \quad V(\bs\varpi_k)\otimes V(\bs\varpi_l), \quad\text{for } k\le l,\text{ and}\quad V(\bs\pi_k)\otimes V(\bs\varpi_l) \quad\text{for all }k,l.
	\end{equation*}
	Then, $V(\bs{\pi})\otimes V(\bs{\varpi})$ is highest-$\ell$-weight and, moreover, if all the above tensor products are simple, then so is $V(\bs{\pi})\otimes V(\bs{\varpi})$.\hfil\qed
\end{cor}

It is well-known (see \cite{ohscr:simptens} and references therein) that, given $(i,r),(j,s)\in I\times  \mathbb Z_{>0}$ and $a,b\in\mathbb Z$, there exists a finite set $\mathscr R_{i,j}^{r,s} \subseteq \mathbb Z_{>0}$ such that
\begin{equation}\label{e:defredset}
	V(\bs\omega_{i,a,r})\otimes V(\bs\omega_{j,b,s}) \text{ is reducible}\qquad\Leftrightarrow\qquad |d_ia-d_jb| \in \mathscr R_{i,j}^{r,s}.
\end{equation}
Moreover, in that case,
\begin{equation}\label{e:krhwtp}
	V(\bs\omega_{i,a,r})\otimes V(\bs\omega_{j,b,s}) \text{ is  highest-$\ell$-weight}\qquad\Leftrightarrow\qquad d_ia>d_jb.
\end{equation}
If $r=s=1$, we simplify notation and write $\mathscr R_{i,j}$ for $\mathscr R_{i,j}^{1,1}$. 
In particular, \eqref{e:defredset} implies $V(\bs\omega_{i,a,r})$ is real for all $i\in I,a,r\in\mathbb Z,r>0$. 
It follows from \Cref{p:sinter} and \eqref{e:dualtp} that
\begin{equation}\label{e:redsim}
	\mathscr R_{j,i}^{s,r} = \mathscr R_{i,j}^{r,s} = \mathscr R_{i^*,j^*}^{r,s}.
\end{equation}
If $\lie g$ is of type $A$  and $i,j\in I, r,s\in\mathbb Z_{>0}$, the following is well-known: 
\begin{equation}\label{e:typeAR}
	\mathscr R_{i,j}^{r,s} = \{r+s+d(i,j)-2p: - d([i,j],\partial I)\le p<\min\{r,s\}  \}. 
\end{equation}
Here, $[i,j]$ denotes the subdiagram consisting of the nodes between $i,j\in I$ (including $i$ and $j$), $d(i,j) = \#[i,j]-1$ is the distance between the nodes $i$ and $j$, $\partial J$ is the set of monovalent nodes of a connected subdiagram $J$, and $d(J,K) = \min\{d(j,k): j\in J,k\in K\}$ for $J,K\subseteq I$. Since $\{i\}\subseteq I$ is a subdiagram of type $A_1$, it  follows  that 
\begin{equation}\label{e:sl2inRij}
    d_i\,\mathscr R^{r,s}\subseteq \mathscr R_{i,i}^{r,s}.
\end{equation} 
In fact, we have
\begin{equation}\label{e:redsl2sd}
    V(\bs\omega_{i,a,r})_{\{i\}}\otimes V(\bs\omega_{i,b,s})_{\{i\}} \ \ \text{is reducible}\ \ \Leftrightarrow\ \ |a-b|\in \,\mathscr R^{r,s}. 
\end{equation}
	
The following consequence of \eqref{e:krhwtp} and \Cref{c:hlwquot} will often be used.
\begin{cor}\label{c:hlwpqf}
	Let $\bs{\pi},\bs\varpi\in\mathcal{P}^+$ and suppose  there exist pseudo $q$-factorizations
	\begin{equation*}
		\bs\pi = \prod_{k=1}^m \bs\pi_k, \quad  \bs\varpi = \prod_{k=1}^{m'} \bs\varpi_k, 
	\end{equation*}
	such that $V(\bs\pi_k)\otimes V(\bs\varpi_l)$ is highest-$\ell$-weight for all $k,l$.
	Then, $V(\bs{\pi})\otimes V(\bs{\varpi})$ is highest-$\ell$-weight.
\end{cor}
	
\begin{proof}
	Write $\bs\pi_k = \bs\omega_{i_k,a_k,r_k}$. Up to reordering, we can assume $a_k\ge a_l$ for all $k\le l$. It then follows from  \eqref{e:krhwtp} that $V(\bs\pi_{k})\otimes V(\bs\pi_l)$ is highest-$\ell$-weight for all $k\le l$. Similarly, we can assume  $V(\bs\varpi_{k})\otimes V(\bs\varpi_l)$ is highest-$\ell$-weight for all $k\le l$. A direct application of \Cref{c:hlwquot} completes the proof.
\end{proof}

The following will play a crucial role in the proof of some of our main results.
	
\begin{lem}\label{l:3aline}
	If $\lie g$ is of type $A$ and $\bs\omega,\bs\varpi\in\mathcal{KR}$, the module $V(\bs\omega)\otimes V(\bs\omega\bs\varpi)$ is simple.\qed
\end{lem}

\begin{proof}
	The result is clear if $V(\bs\omega)\otimes V(\bs\varpi)$ is simple. If this is not the case and  $\bs\omega,\bs\varpi$ are the $q$-factors of $\bs\pi:=\bs\omega\bs\varpi$, this is exactly the statement of  \cite[Corollary 2.4.9]{ms:3tree} (which is a particular case of \Cref{t:3lineprime} reviewed below). Otherwise, 
	the well-known combinatorics of $q$-strings implies that $\bs\pi$ has at most two $q$-factors: the one with higher degree is divisible by  $\bs\omega$ and the other divides $\bs\omega$. If we denote them by $\bs\omega_1$ and $\bs\omega_2$ (allowing the possibility that $\bs\omega_2=\bs 1$), it follows
	that $V(\bs\pi)\cong V(\bs\omega_1)\otimes V(\bs\omega_2)$ and $V(\bs\omega)\otimes V(\bs\omega_j)$ is simple for $j\in\{1,2\}$. Hence, $V(\bs\pi)\otimes V(\bs\omega)$ is also simple.
\end{proof}
	
\subsection{KKOP Invariants}\label{ss:KKOP}
A non-negative integer  $\lie d(V, W)$ was defined in \cite{kkop:mc1}  for each pair of simple modules $V$ and $W$. It follows directly from the definition that 
\begin{equation*}
	\lie d(V, W) = \lie d(W,V).
\end{equation*}
We summarize important properties of $\lie d$ which we will use several times for checking the simplicity of certain tensor products. We shall refer to the number $\lie d(V,W)$ as the KKOP invariant of the pair $V,W$. 
	
\begin{prop}\label{kkop} 
	Let $V_1$, $V_2$, and $W$ be simple modules and set $V=V_1\otimes V_2$. Then for any subquotient $S$ of $V$ we have 
	\begin{equation}\label{e:kkopineq}
		\lie d(S, W)\le \lie d(V_1,W) +\lie d(V_2, W).
	\end{equation}
	Moreover, if either $V_1$ or $V_2$ is real the following hold.
	\begin{enumerate}[(a)]
		\item $V$ is simple if and only if  $\lie d(V_1, V_2)=0$.
		\item If $\lie d(V_1, V_2)=1$ then $V$  has length two.
		\item If $V_1$ and $V_2$ are real  and $\lie d(V_1,V_2)\leq 1$ then ${\rm soc }(V)$ and ${\rm hd}(V)$ are real.		\qed
	\end{enumerate}
\end{prop}

Equation \eqref{e:kkopineq} follows from \cite[ Proposition 4.2 and Lemma 3.10]{kkop:mc1}, part (a) follows from \cite[Corollary 3.17]{kkop:mc1} and part (b) from   \cite[Proposition 4.7]{kkop:mc1} and (c) from \cite[Lemma 2.27 and Lemma 2.28]{kkop:pbw}.

\begin{cor}\label{c:kkop}
	Let $V_1,V_2$, and $V$ be as in \Cref{kkop}. If $V_1$ is real and $\lie d (V_1,V_2)=1$, then
	\begin{equation*}
		{\rm soc }(V)\otimes V_1 \quad\text{and}\quad {\rm hd}(V)\otimes V_1 \quad\text{are simple}.
	\end{equation*}
\end{cor}

\begin{proof}
	Since $\lie d(V_1,V_2)=1>0$ it follows from \cite[Corollary 4.12]{kkop:mc1} that 
	$$\lie d({\rm soc}(V),V_1)=0=\lie d({\rm hd}(V),V_1).$$
	The result is now immediate from Proposition \ref{kkop}(a) since $V_1$ is real.
\end{proof}

\begin{lem}[{\cite[Lemma 6.2.1]{naoi:Tsys}}]\label{l:naoi} 
	If $\lie g$ is of type $A, k\ge 2$, and $\bs\pi = \bs\omega_{i_1,a_1}\cdots\bs\omega_{i_k,a_k}$ is such that 
	\begin{equation*}
		a_{s+1}-a_s\in \mathscr R_{i_s,i_{s+1}} \quad\text{for all}\quad 1\leq s<k,
	\end{equation*}
	then $\lie d(V(\bs\omega_{i_1,a_1}), V(\bs\pi\bs\omega_{i_1,a_1}^{-1}))=1 = \lie d (V(\bs\omega_{i_k,a_k}),V(\bs\pi\bs\omega_{i_k,a_k}^{-1}))$.\qed
\end{lem}

\subsection{Digraph Notation and Terminology}\label{ss:digraphs} By a digraph (or simply a graph) $G$ we shall mean a pair $(\mathcal V_G,\mathcal A_G)$, where $\mathcal V_G$ is a finite set and $\mathcal A_G\subseteq \mathcal V_G\times\mathcal V_G$. When no confusion arises, we shall simply write $\mathcal V$ and $\mathcal A$. The elements of $\mathcal V$ are called the vertices of $G$ and the ones in $\mathcal A$ are called the arrows. If $a=(v,w)\in\mathcal A$, then $v$ is said to be the tail of $a$, denoted $t(a)$, and $w$ is the head of $a$, denoted $h(a)$. 
	
A pair $(\mathcal V', \mathcal A')$ is said to be a subgraph of $G$ if $\mathcal V'\subseteq \mathcal V$ and $\mathcal A' = \{a\in\mathcal A: t(a),h(a)\in\mathcal V'\}$. If $\mathcal U\subseteq \mathcal V$, we denote by $G_{\mathcal U}$ the subgraph of $G$ whose vertex set is $\mathcal U$. If $H=G_{\mathcal U}$, we use the notation $H^c$ as well as $G\setminus H$ for the subgraph $G_{\mathcal V\setminus\mathcal U}$. If $\mathcal U$ has a single element $u$, we may simplify notation and write $G_u$ for $G_{\mathcal U}$ and $G\setminus u$ for $G\setminus G_{\mathcal U}$. We shall write $H\triangleleft G$ to say $H$ is a subgraph of $G$.
	
A pair of subgraphs of the form $(H,H^c)$ is often called a cut of $G$ and the set
\begin{equation*}
	\mathcal A\setminus  (\mathcal A_H\cup\mathcal A_{H^c})
\end{equation*}
is called the associated cut set. We shall say an arrow in the cut set links the subgraphs $H$ and $H^c$. If $H$ and $H'$ are subgraphs of $G$ with  vertex sets $\mathcal U$ and $\mathcal U'$, we define 
\begin{equation}
	H\otimes H'=G_{\mathcal U\,\cup\,\mathcal U'}.
\end{equation}
If $\mathcal U$ and $\mathcal U'$ are disjoint, then $(H,H')$ is a cut of $H\otimes H'$. Thus, for disjoint subgraphs $H$ and $H'$, we shall say an arrow links $H$ and $H'$ in $G$ if it is an element of the cut-set of $H\otimes H'$ associated to $(H,H')$.    The notion of tensor product of any finite family of subgraphs is defined in the obvious way and it is clearly associative and commutative.

More generally, let $\mathcal G=G_1,\dots, G_l$ be a sequence of subgraphs of a graph $G$ and let $\mathcal V_k$ and $\mathcal A_k$ be the corresponding vertex and arrow sets. We shall say $\mathcal G$ is a cut of length $l$ if
\begin{equation*}
	\mathcal V = \bigcup_{k=1}^l\mathcal V_k \quad\text{and}\quad \mathcal V_k\cap \mathcal V_m = \emptyset\ \ \text{for all}\ \ 1\le k<m\le l.
\end{equation*}
Thus, if $l=1$ we have $G_1=G$ and, if $l=2$, we recover the usual definition of cut. We shall simply say ``$\mathcal G$ is a multicut'' if $l\ge 2$ and we are not interested in being explicit about the length. The corresponding cut-set is defined as
\begin{equation}\label{e:cutsetgen}
	\mathcal A_{\mathcal G} = \mathcal A \setminus \bigcup_{k=1}^l \mathcal A_k.
\end{equation}
If $\mathcal G$ is a cut of length $l$, set
\begin{equation}
	\bar G_k = G_{k+1}\otimes\cdots\otimes G_l  \quad\text{for}\quad 0\le k\le l. 
\end{equation}
Note 
\begin{equation}\label{e:whyquochain}
	G_{k}\triangleleft \bar G_{k-1}, \quad  \bar G_k = \bar G_{k-1}\setminus G_k \quad\text{for all}\quad 0< k\le l,
\end{equation}
and the sequence $\bar G_0,\dots,\bar G_l$ is a proper descending chain of subgraphs:
\begin{equation}\label{e:defquochain}
	\emptyset=\bar G_l \triangleleft\bar G_{l-1}\triangleleft\cdots\triangleleft \bar G_1\triangleleft\bar G_0 = G.
\end{equation}
We shall refer to this chain as the quochain associated to $\mathcal G$. By abuse of language, we shall often refer to $\mathcal G$ as a quochain as well. Let us explain the choice for the term ``quochain''. One interpretation of \eqref{e:whyquochain} is that $\bar G_k$ is the ``quotient'' of $\bar G_{k-1}$ by the (tensor) factor $G_k$. In other words, a quochain is a chain of quotients. The notion of quochains will play a prominent role in the main results of this paper.

Suppose $P$ is a ``graphical property'', i.e., a property assignable to subgraphs of a graph, such as being connected. We shall say the multicut $\mathcal G$ determines a $P$-quochain (or that $\mathcal G$ is a $P$-quochain by abuse of language) if $G_k$ has the property $P$ when regarded as a subgraph of $\bar G_{k-1}$ for all $1\le k\le l$. We remark that $P$ may not be ``intrinsic''. For instance, one may work with graphs equipped with extra structure (as it will be the case here) and $P$ may be related to the extra structure. We will use this notion with a combination of intrinsic and extrinsic properties $P$.

We now fix a graph $G$ and set some terminology related to vertices of  $G$.
Given $v\in \mathcal V$, set
\begin{equation}\label{e:adjnot}
	\hat A_v = \{w\in\mathcal V:(v,w)\in\mathcal A\}, \quad \check A_v = \{w\in\mathcal V:(w,v)\in\mathcal A\}, \quad \bar A_v = \hat A_v\cup\check A_v.
\end{equation}
A vertex in $\bar A_v$ is said to be adjacent to $v$. The valence of $v$ is defined as $\#\bar A_v$. A vertex $v$ is said to be a source if $\check A_v=\emptyset$ and it is a sink if $\hat A_v=\emptyset$.  Note that valence-1 vertices are necessarily either sinks or sources. The adjacency subgraph of $v$ is defined as 
\begin{equation}
	\Adj(v) = G_{\{v\}\cup\bar A_v}.
\end{equation}
When $v$ is a vertex of a subgraph $H$, we shall write $\Adj_G(v)$ and $\Adj_H(v)$ for the corresponding adjacency subgraphs. Evidently,
\begin{equation*}
	\Adj_H(v)\subseteq \Adj_G(v),
\end{equation*}
but we may have proper inclusions.
	
A loop is a an arrow $a$ such that $t(a)=h(a)$. A path $\rho$ in $G$ is a sequence of vertices $v_1,\dots,v_l$ such that $\{(v_j,v_{j+1}), (v_{j+1},v_j)\}\cap \mathcal A\ne\emptyset$ for all $1\le j< l$. The path is said to be oriented if  either $(v_j,v_{j+1})\in \mathcal A$ for all $1\le j< l$, or $(v_{j+1},v_j)\in\mathcal A$ for all $1\le j< l$. It is said to contain a cycle if $v_j = v_{j'}$ for some $1\le j\neq j'\le l$. The graph $G$ is said to be connected if there is a path linking any two of its vertices. 

If there are no loops nor oriented cycles in $G$, the set $\mathcal A$ induces a partial order on $\mathcal V$ by transitive extension of the relation $h(a)<t(a)$ for all $a\in\mathcal A$. Henceforth, we assume $G$ has this property. In this case, $v$ is a source if and only if $v$ is a maximal element for this partial order. Moreover, $G$ has at least one source and one sink. We shall say $G$ is totally ordered if this order is linear  (such graphs are also called traceable in the literature). A maximal totally ordered subgraph is a subgraph whose vertex set is not properly contained in the vertex set of another totally ordered subgraph.
We shall say $H$ is an extremal subgraph if 
\begin{equation}\label{e:extsg}
	\begin{aligned}
		\text{either}\quad  & t(a)\in\mathcal V_H \text{ for all } a\in\mathcal A\setminus(\mathcal A_H\cup\mathcal A_{H^c})\\ \quad\text{or}\quad & h(a)\in\mathcal V_H \text{ for all } a\in\mathcal A\setminus(\mathcal A_H\cup\mathcal A_{H^c}).
	\end{aligned}
\end{equation}  
If the first option above occurs, we shall say $H$ is top subgraph. Otherwise, we will say it is a bottom subgraph. Evidently, $H$ is a bottom subgraph if and only if $H^c$ is a top subgraph.

\section{A Reality Test via Graphs}\label{s:realbygraph}
			
\subsection{Pseudo $q$-Factorization Graphs}\label{ss:qfgraphs}
The notion of pseudo $q$-factorization graphs was introduced in \cite{ms:to}. We shall review it starting with a rephrasing of the original definition. Recall that $\mathcal {KR}$ denotes the set of KR type Drinfeld polynomials as well as \eqref{e:krhwtp}.
			
Let $G=(\mathcal V,\mathcal A)$ be a digraph. A pseudo $q$-factorization map over $G$ is a map $\mathcal F:\mathcal V\to\mathcal{KR}$ such that
	\begin{equation}
		\mathcal F(v) = \bs\omega_{i,a,r} \ \ \text{and}\ \ \mathcal F(w) = \bs\omega_{j,b,s} \quad\Rightarrow\quad \Big[ (v,w)\in\mathcal A \quad\Leftrightarrow\quad d_ia-d_jb \in\mathscr R_{i,j}^{r,s} \Big].
	\end{equation}
	A pseudo $q$-factorization graph is a digraph equipped with a pseudo $q$-factorization map.  If $H=G_{\mathcal U}$ is a subgraph, then $\mathcal F|_{\mathcal U}$ turns $H$ into a pseudo $q$-factorization graph. Specifying the map $\mathcal F$ is equivalent to specifying three maps: $c:\mathcal V\to I,\lambda:\mathcal V\to\mathbb Z_{>0},\epsilon:\mathcal A\to \mathbb Z_{>0}$, up to a uniform shift of all the centers of the $q_i$-strings for each connected component of $G$. These maps are determined by:
	\begin{equation}\label{e:defcw}
		\mathcal F(v) = \bs\omega_{i,a,r},  \quad\Rightarrow\quad c(v)= i, \ \ \lambda(v) = r,
	\end{equation}
	and
	\begin{equation}\label{e:defexp}
		\mathcal F(v) = \bs\omega_{i,a,r}, \ \ \mathcal F(w) = \bs\omega_{j,b,s}, \ \ (v,w)\in\mathcal A \quad\Rightarrow\quad \epsilon(v,w) = d_ia-d_jb.
	\end{equation}
	It is not difficult to check using \eqref{e:defexp} that $\epsilon$ is compatible with paths in the sense of \cite[Eq. (2.5.1)]{ms:to}. In particular, if a pseudo $q$-factorization map on $G$ exists, $G$ does not contain loops nor oriented cycles. 

	Conversely, suppose $c:\mathcal V\to I,\lambda:\mathcal V\to\mathbb Z_{>0},\epsilon:\mathcal A\to \mathbb Z_{>0}$ are given and satisfy
	\begin{equation*}
		(v,w)\in\mathcal A \quad\Leftrightarrow\quad \epsilon(v,w) \in\mathscr R_{c(v),c(w)}^{\lambda(v),\lambda(w)}.
	\end{equation*}
	Suppose further that $\epsilon$ is compatible with paths in the sense of \cite[Eq. (2.5.1)]{ms:to} and that $G$ is connected. Then, given $v\in\mathcal V$ and $a\in\mathbb Z$, there exists a unique pseudo $q$-factorization map $\mathcal F$ over $G$ which induces $c,\lambda,\epsilon$ in the sense of \eqref{e:defcw} and \eqref{e:defexp}, and such that $\mathcal F(v)=\bs\omega_{c(v),a,\lambda(v)}$. Indeed, if $(v,w)\in\mathcal A$, then we must have (cf. \eqref{e:defexp})
	\begin{equation*}
		\mathcal F(w) = \bs\omega_{c(w),\frac{1}{d_{c(w)}}(d_{c(v)}a-\epsilon(v,w)),\lambda(w)}.
	\end{equation*}
	If it is $(w,v)\in\mathcal A$ instead, then we must replace the minus sign above by a plus sign. The assumption that $G$ is connected implies $\mathcal F$ is completely determined by this rule. The maps $c,\lambda,\epsilon$ are called the color, the weight, and the exponent components of $\mathcal F$.

We shall say that a pseudo $q$-factorization map $\mathcal F$ over $G$ is fundamental if $\lambda(v)=1$ for all $v\in\mathcal V$.  The corresponding pseudo $q$-factorization graph will then be referred to as a fundamental factorization graph. Recall \eqref{e:defqfact} and \eqref{e:sl2inRij}. We shall say a pseudo $q$-factorization map $\mathcal F$ is a $q$-factorization map if 
\begin{equation}
	c(v) = c(w) \ \ \text{and}\ \ (v,w)\in\mathcal A \quad\Rightarrow\quad \epsilon(v,w) \notin d_{c(v)}\,\mathscr R^{\lambda(v),\lambda(w)}.
\end{equation}
The corresponding pseudo $q$-factorization graph will then be referred to as a $q$-factorization graph.
			
Given a pseudo $q$-factorization map $\mathcal F$ over $G$, define
\begin{equation}\label{e:polyofgraph}
	\bs\pi_{\mathcal F} = \prod_{v\in\mathcal V} \mathcal F(v)\in\mathcal P^+.
\end{equation}
Then, the right-hand-side is a pseudo $q$-factorization of $\bs\pi_{\mathcal F}$ by construction. By abuse of notation, we shall identify $v\in\mathcal V$ with $\mathcal F(v)$, so we can shorten the above to $\bs\pi_{\mathcal F} =  \prod_{v\in\mathcal V} v$. Moreover, we shall abuse of language and simply say ``$G$ is a pseudo $q$-factorization graph'' with no mention to the structure data $(\mathcal V,\mathcal A,\mathcal F)$ and then write $\bs\pi_G$ instead of $\bs\pi_{\mathcal F}$. Despite the lack of accuracy, this should be most often beneficial to the writing of the text. We shall also say $G$ is a pseudo $q$-factorization graph over $\bs\pi$ if $\bs\pi_G=\bs\pi$. If $H=G_{\mathcal U}$ is a subgraph, we set
\begin{equation}\label{e:polyofsubgraph}
	\bs\pi_{\mathcal U} = \bs\pi_H = \bs\pi_{\mathcal F|_{\mathcal U}}.
\end{equation}
In particular,
\begin{equation*}
	\bs\pi_G = \bs\pi_H\, \bs\pi_{H^c}.
\end{equation*}

Conversely, given any map $\mathcal F:\mathcal V\to\mathcal{KR}$ defined on a nonempty finite set $\mathcal V$,  we can construct a pseudo $q$-factorization graph having $\mathcal V$ as vertex set and $\mathcal F$ as its pseudo $q$-factorization map by defining $\mathcal A$ by the requirement:
\begin{equation*}
	(v,w)\in\mathcal A \quad\Leftrightarrow\quad V(\mathcal F(v))\otimes V(\mathcal F(w)) \ \ \text{is reducible and highest-$\ell$-weight.}
\end{equation*}
We denote this graph by $G(\mathcal F)$. If $\mathcal F$ is an actual $q$-factorization map over $G(\mathcal F)$ and $\bs\pi=\bs\pi_{\mathcal F}$, we also use the notation $G(\bs\pi)$ and call it the $q$-factorization graph of $\bs\pi$.  If $G(\mathcal F)$ is a fundamental factorization graph, we also use the notation $G_f(\bs\pi)$ and call it the fundamental (factorization) graph of $\bs\pi$. 
			
Given pseudo $q$-factorization graphs $G$ and $G'$ over $\bs\pi$ and $\bs\pi'$, respectively, we denote by $G\otimes  G'$ the unique pseudo $q$-factorization graph over $\bs\pi\bs\pi'$ whose vertex set is $\mathcal V_G\, \dot{\cup}\, \mathcal V_{G'}$. In particular, $(G,G')$ is a cut of $G\otimes G'$, so the notation agrees with that defined in \Cref{ss:digraphs}. Note however that $G\otimes G'=G'\otimes G$, but this has no relation with whether $V(\bs\pi_G)\otimes V(\bs\pi_{G'})\cong V(\bs\pi_{G'})\otimes V(\bs\pi_{G})$ or not.
In the case that $G'$ has a single vertex $v$, we shall also use the notation $G\otimes v$ for $G\otimes G'$. We shall say that $G$ and $G'$ are dissociate pseudo $q$-factorization graphs if the corresponding cut set is empty, i.e, there are no arrows in $G\otimes G'$ linking $G$ and $G'$. In that case, \Cref{c:hlwquot} implies
\begin{equation}\label{e:dissosimple}
	V(\bs\pi)\otimes V(\bs\pi') \quad\text{is simple.}
\end{equation}
This implies that the study of primality and reality of modules can be reduced to considering connected pseudo $q$-factorization graphs. If $v\in\mathcal V_G$, we shall say $G$ and $G'$ are linked through $v$ if there exists an arrow $a$ of $G\otimes G'$ linking $G$ and $G'$ such that
\begin{equation}\label{e:linkedparts}
	v\in\{h(a),t(a)\}.
\end{equation}
			
Let $P$ be a representation theoretic property, i.e., a property assignable to modules such as being prime or real. We shall interpret $P$ as a (extrinsic) graphical property as follows. Suppose $H$ is a subgraph of a pseudo $q$-factorization graph $G$. We shall say $H$ satisfies $P$ if $V(\pi_H)$ satisfies $P$. For instance, we shall say $H$ is real if $V(\pi_H)$ is real. 
Suppose $\mathcal G=G_1,\dots, G_l$ is a multicut of $G$. Given a graphical property $P$, we shall say $\mathcal G$ has the property $P$ if $G_k$ has property $P$ for all $k$ (when regarded as subgraphs of $G$). For instance, if both $H$ and $H^c$ are real subgraphs of $G$, we shall say the associated cut is real and that $H$ determines a real cut of $G$. 

Let $\kappa^*$ be the group automorphism of $\mathcal P$ determined by $\bs\omega_{i,a}\mapsto \bs\omega_{i,-a}$. 
Given a pseudo $q$-factorization graph $G$, it follows from \eqref{e:cart-} and the comments preceding it that
\begin{equation}\label{e:arrowd}
	{}^*(V(\bs\pi_G)^\kappa) \cong V(\kappa^*(\bs\pi_G)). 
\end{equation}
Let $G^{\rm opp}$ be the digraph $(\mathcal V,\mathcal A^{\rm opp})$, where 
\begin{equation*}
	(v,w)\in \mathcal A^{\rm opp} \quad\Leftrightarrow\quad (w,v)\in\mathcal A.
\end{equation*}
If $\mathcal F$ is the pseudo $q$-factorization map over $G$, one easily checks using \eqref{e:redsim} that $\mathcal F^- = \kappa^*\circ\mathcal F$ is a pseudo $q$-factorization map over $G^{\rm opp}$ and
\begin{equation*}
	\bs\pi_{G^{\rm opp}} = \kappa^*(\bs\pi_G).
\end{equation*}
Evidently, $V(\bs\pi_G)$ is real (prime, etc) if and only if so is $V(\bs\pi_{G^{\rm opp}})$. We shall refer to the map $G\mapsto G^{\rm opp}$ as arrow duality. 
			
\subsection{On the Prime Factorization of 3-Vertex Non Totally Ordered Graphs}\label{ss:3v}		
We recall in this section one of the main results of \cite{ms:3tree} which will be used in a couple of examples below. Thus, let $G$ be a $q$-factorization graph with $3$ vertices. If $G$ is not totally ordered, then it must be of the form
\begin{equation}\label{e:3valg}
	\begin{tikzcd}
		\stackrel{r_1}{i_1} & \arrow[swap,l,"m_1"]  \stackrel{r}{i} \arrow[r,"m_2"] & \stackrel{r_2}{i_2} 
	\end{tikzcd} \qquad\text{or}\qquad 
	\begin{tikzcd}
		\stackrel{r_1}{i_1} \arrow[r,"m_1"] &   \stackrel{r}{i}  & \arrow[swap,l,"m_2"] \stackrel{r_2}{i_2} 
	\end{tikzcd}
\end{equation}
Here, $i$ and $i_j$ are the corresponding coloring of the vertices, $r$ and $r_j$ are the corresponding weights, while $m_j$ are the corresponding exponents. Assume henceforth that $\lie g$ is of type $A$.

Given $i,j\in I, r,s\in\mathbb Z_{>0}$, and a connected subdiagram $J$ containing $\{i,j\}$, let $\mathscr R_{i,j,J}^{r,s}$ be determined by 
\begin{equation*}
	V((\bs\omega_{i,a,r})_J)\otimes V((\bs\omega_{j,b,s})_J)  \text{ is reducible}\qquad\Leftrightarrow\qquad a-b = m \text{ with } |m|\in \mathscr R_{i,j,J}^{r,s}.
\end{equation*}
It follows from \eqref{e:sJs} that $\mathscr R_{i,j,J}^{r,s}\subseteq \mathscr R_{i,j,K}^{r,s}$ if $J\subseteq K$.

\begin{thm}[{\cite[Theorem 2.4.6]{ms:3tree}}]\label{t:3lineprime}
	Assume $\lie g$ is of type $A$ and let $G=G(\bs\pi)$ be as in \eqref{e:3valg}. For $j=1,2$, let also $I_j\subseteq I$ be the minimal connected subdiagram containing $\{i,i_j\}$ such that $m_j\in\mathscr R_{i,i_j,I_j}^{r,r_j}$ and let $j'$ be such that $\{j,j'\}=\{1,2\}$. Then, $G$ is not prime if and only if there exists $j\in\{1,2\}$ such that
	\begin{equation*}\label{2gencondm}
		i_{j'}\in I_j, \qquad m_{j'}\in{\mathscr{R}_{i,i_{j'},I_j}^{r,r_{j'}}}, \qquad m_{j'}-m_j+\check h_{I_j}\in\mathscr R_{w_0^{I_j}(i_j),i_{j'},I_j}^{r_j,r_{j'}},
	\end{equation*}
	and
	\begin{equation}\label{2exracondm}
		m_j+r_j \le m_{j'} + r_{j'} +d(i_1,i_2).
	\end{equation}
	In that case, $V(\bs\pi)\cong V(\bs\omega)\otimes V(\bs\pi\bs\omega^{-1})$, where $\bs\omega$ is the $q$-factor corresponding to such $j$.\qed
\end{thm}

\subsection{Strongly Real Modules}\label{ss:realtest}

The following theorem motivates the main definitions of this paper: that of reality determining subgraphs and that of strongly real modules. These concepts will allow us to describe certain families of real modules.

\begin{thm}\label{t:isreal}
	Suppose  $\bs\pi_1,\bs\pi_2\in \mathcal P^+$ satisfy the following:
	\begin{enumerate}[(i)]
		\item $V(\bs\pi_j)$ is real for $j=1,2$;
		\item $V(\bs\pi_1)\otimes V(\bs\pi_2)$ is highest $\ell$-weight;
		\item $V(\bs\pi_1\bs\pi_2)\otimes V(\bs\pi_j)$ is simple for some $j=1,2$.
	\end{enumerate}
	Then, $V(\bs\pi_1\bs\pi_2)$ is real. 
\end{thm}
			
\begin{proof}
	We shall write the details assuming (iii) holds for $j=1$. The other case is proved similarly. To shorten notation, set $\bs\pi=\bs\pi_1\bs\pi_2$. 
				
	Since $V(\bs\pi_1)$ is real, $V(\bs\pi_1^k)\cong V(\bs\pi_1)^{\otimes k}$. It then follows from (ii) and \Cref{t:cyc} that 
	\begin{equation*}
		V(\bs\pi_1^k)\otimes V(\bs\pi_2) \quad\text{is highest-$\ell$-weight for all $k>0$,}
	\end{equation*}	
	and, therefore,  there exists an epimorphism 
	\begin{equation*}
		V(\bs\pi_1^k)\otimes V(\bs\pi_2)\to V(\bs\pi_2\bs\pi_1^k).
	\end{equation*}
	Hence,  there also exists an epimorphism	
	\begin{equation*}
		V(\bs\pi_1^2)\otimes V(\bs\pi_2)\otimes  V(\bs\pi_2) \to 	V(\bs\pi_2\bs\pi_1^2)\otimes V(\bs\pi_2).
	\end{equation*}
	Since all the tensor factors on the left-hand-side above are real, \Cref{t:cyc} implies
	\begin{equation}\label{e:realth}
		V(\bs\pi_2\bs\pi_1^2)\otimes V(\bs\pi_2)\quad\text{is highest-$\ell$-weight.}
	\end{equation}
	Finally, we have an epimorphism
	\begin{equation*}
		V(\bs\pi_2\bs\pi_1^2)\otimes V(\bs\pi_2)\stackrel{(iii)}{\cong} V(\bs\pi)\otimes V(\bs\pi_1)\otimes V(\bs\pi_2) \to V(\bs\pi)\otimes V(\bs\pi)
	\end{equation*}
	and \eqref{e:realth} then implies $V(\bs\pi)\otimes V(\bs\pi)$ is highest-$\ell$-weight. Reality then follows from \Cref{p:vnvstar}.	
\end{proof}

Recall the paragraph after \eqref{e:linkedparts}, especially the definition of real cut. 

\begin{defn}\label{d:rds}
	A subgraph $H$ of a pseudo $q$-factorization graph $G$ will be said a reality determining subgraph (rds for short) if, either $\#\mathcal V_H=\#\mathcal V_G=1$, or $H$ is a proper nonempty subgraph satisfying the following:
	\begin{enumerate}[(i)]
		\item $H$ determines a real cut;
		\item Either $V(\bs\pi_H)\otimes V(\bs\pi_{H^c})$ or $V(\bs\pi_{H^c})\otimes V(\bs\pi_{H})$ is highest-$\ell$-weight;
		\item $V(\bs\pi_G)\otimes V(\bs\pi_H)$ is simple.\endd
	\end{enumerate}
\end{defn}	
		
It follows from \Cref{t:isreal} that $G$ is real, i.e., $V(\bs\pi_G)$ is real, if $G$ contains an rds. Condition (i) will usually be checked from previous examples (recursive procedures, for instance). The easiest situation for checking condition (ii) is that of extremal subgraphs in light of \Cref{l:topbothw} below. Thus, in most examples we consider, the hard part of checking if a subgraph is an rds is condition (iii).
			
\begin{lem}\label{l:topbothw}
	Let $H$ and $H'$ be subgraphs with disjoint vertex-sets of a pseudo $q$-factorization graph $G$ such that $H$ is a top subgraph of $H\otimes H'$. Then, $V(\bs\pi_H)\otimes V(\bs\pi_{H'})$ is highest-$\ell$-weight.
\end{lem}

\begin{proof}
    Given vertices $v\in\mathcal V_H$ and $v'\in\mathcal V_{H'}$, the assumption that $H$ is a top subgraph (recall the definition after \eqref{e:extsg}) implies that either there is no arrow between $v$ and $v'$ or $(v,v')$ is a arrow. By definition of pseudo $q$-factorization, in the former case $V(\mathcal F(v))\otimes V(\mathcal F(v'))$ is simple, while it is reducible and highest $\ell$-weight in the latter. In any case, it is highest-$\ell$-weight and the claim follows from \Cref{c:hlwpqf}.
\end{proof}
			
\begin{cor}\label{c:isreal}
	Let $G$ be a pseudo $q$-factorization over $\bs\pi$ and suppose there exists an extremal subgraph $H\triangleleft G$ which determines a real cut and satisfies
	\begin{equation}\label{e:realts}
		V(\bs\pi)\otimes V(\bs\pi_H) \quad\text{is simple.}
	\end{equation}
	Then, $H$ is an rds. 
\end{cor}
			
\begin{proof}
	If $G$ is a singleton, there is nothing to be done. Otherwise, let $K=H$ if $H$ is a top subgraph and $K=H^c$ otherwise. Set $\bs\pi_1=\bs\pi_K$ and $\bs\pi_2=\bs\pi_{K^c}$. Then, \Cref{l:topbothw} implies (ii)  is satisfied, while (i) follows from the assumption that $H$ determines a real cut and (iii) follows from \eqref{e:realts}. 
\end{proof}

Recall the definition of $P$-quochains given soon after \eqref{e:defquochain}. Thus, to say $G$ admits an rds-quochain (of length $l$) means that there exists a multicut $\mathcal G=G_1,\dots,G_l$ of $G$ such that $G_k$ is an rds in $\bar G_{k-1}$ for all $1\le k\le l$. In that case, $G_l$ is necessarily a singleton.

\begin{defn}\label{d:stronglyreal}
	A pseudo $q$-factorization graph $G$ is said to be strongly real if it admits an rds-quochain for which all the rds are singletons. \endd
\end{defn}

Equivalently, $G$ is strongly real if it admits an rds-quochain whose length is $\#\mathcal V_G$.

\begin{defn}\label{d:weakrdsq}
	A weak rds-quochain of length $l$ for $G$ is a multicut   $\mathcal G=G_1,\dots,G_l$  of $G$ such that  $G_k$ is an rds in $\bar G_{k-1}$ for all $1\le k< l$ and $G_l$ is real.\endd
\end{defn}

Note that, if $\mathcal G$ is a weak rds-quochain of length $l>1$, the assumption that $G_l$ is real in the above definition follows from the assumption that  $G_{l-1}$ determines a real cut of $G_{l-1}\otimes G_l=\bar G_{l-2}$. In other words, the last assumption of the definition is necessary only for $l=1$. In particular, the real graphs with no rds are exactly those with at least two vertices and whose weak rds-quochains have length one. We make further comments in this direction in \cref{ss:nords}.

We will say the simple module $V(\bs\pi)$ is strongly real if there exists a pseudo $q$-factorization graph $G$ over $\bs\pi$ which is strongly real. Note that if $V(\bs\pi)$ is strongly real, then it is real by the comment after \cref{d:rds} since there exists a pseudo $q$-factorization graph $G$ over $\bs\pi$ containing an rds. As we shall see in the next sections, pseudo $q$-factorization graphs afforded by trees, snake modules, as well as the more general class of snake trees, which we introduce in \Cref{ss:srt}, are examples of strongly real modules. In particular, examples of snake trees which are not snakes nor trees are given in \Cref{ss:srt}. In all of these cases, we will show that we can chose an rds-quochain of singletons such that all the rds in the sequence are extremal subgraphs.

In \Cref{ss:eas} and \Cref{ss:etoas}, we will describe sufficient conditions for an extremal vertex of a pseudo $q$-factorization graph to be an rds. Successive applications of such criteria will then be used to prove the aforementioned snake trees are strongly real.  
			
Let us give an example of a real module whose $q$-factorization graph admits an rds, but it is not strongly real.

\begin{ex}\label{ex:realbutnotstrong}
	Let $\bs\pi=2_01_33_32_6^2$ in type $A_3$. Thus, $G_f(\bs\pi)=G(\bs\pi)$ is the following graph
	
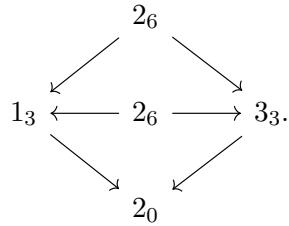
\begin{equation*}
		\begin{tikzcd}
			& 2_6 \ar[rd] \ar[ld] & \\
			1_3 \ar[rd] & 2_6 \ar[r] \ar[l] & 3_3 \ar[ld]. \\
			& 2_0 &
		\end{tikzcd}
	\end{equation*}
	Let $H$ be any of the subgraphs such that 
	\begin{equation*}
		\bs\pi_H=1_32_6.
	\end{equation*} 
	We will check $H$ is an rds and, hence, $V(\bs\pi)$ is real. Since $H$ and $H^c$ are trees, it follows from \cite[Theorem 2.4.8]{ms:3tree} that $H$ determines a real cut. 
	To check that condition (ii) of \Cref{d:rds} holds, we will show
	\begin{equation}\label{e:realbutnotstrongii}
		V(\bs\pi_H)\otimes V(\bs\pi_{H^c}) \quad\text{is highest-$\ell$-weight.}
	\end{equation}
	Then, in order to check that condition (iii) also holds, we will check
	\begin{equation}\label{e:realbutnotstrongiii}
		V(\bs\pi)\otimes V(\bs\pi_H) \quad\text{is highest-$\ell$-weight.}
	\end{equation}
	Note that a combination of \eqref{e:realbutnotstrongii}, \eqref{e:realbutnotstrongiii}, \Cref{t:cyc}, and \Cref{p:vnvstar}, implies that condition (iii) holds. Indeed, \eqref{e:realbutnotstrongii} and \Cref{t:cyc} imply
	\begin{equation*}
		V(\bs\pi_H)\otimes V(\bs\pi_H)\otimes V(\bs\pi_{H^c}) \quad\text{is highest-$\ell$-weight}
	\end{equation*}
	and projects onto $V(\bs\pi_H)\otimes V(\bs\pi)$. Hence, $V(\bs\pi_H)\otimes V(\bs\pi)$ is highest-$\ell$-weight and, hence, alongside with \Cref{p:vnvstar} and \eqref{e:realbutnotstrongiii}, we conclude that $V(\bs\pi_H)\otimes V(\bs\pi)$ is simple, as required by condition (iii).
				
	In order to check \eqref{e:realbutnotstrongii}, begin by noting that  $V(\bs\pi_H)\otimes V(2_6)$ is simple by \cite[Corollary 2.4.9]{ms:3tree}. An application of \Cref{t:cyc}  and \cref{c:hlwquot} then implies 
	\begin{equation*}
		V(\bs\pi_H)\otimes V(2_6)\otimes V(3_3)\otimes V(2_0) \quad\text{is highest-$\ell$-weight}
	\end{equation*}
	and, hence, we have an epimorphism
	\begin{equation*}
		V(\bs\pi_H)\otimes V(2_6)\otimes V(3_3)\otimes V(2_0)\to V(\bs\pi_H)\otimes V(\bs\pi_{H^c}),
	\end{equation*}
	which completes the proof of \eqref{e:realbutnotstrongii}.
				
	For checking \eqref{e:realbutnotstrongiii}, begin by noting that the above argument also shows \begin{equation}\label{e:13263326}
		V(\bs\pi_H)\otimes V(2_63_3) \quad\text{is highest-$\ell$-weight.}
	\end{equation}
	In particular, another application of \Cref{t:cyc} implies
	\begin{equation*}
		V(2_61_3)\otimes V(2_63_3)\otimes V(2_0)\otimes V(2_6)\quad\textrm{and}\quad V(2_6^23_3)\otimes V(1_32_0)\otimes V(1_3)
	\end{equation*}
	are both highest-$\ell$-weight and we have epimorphisms
	\begin{gather*}
		V(2_61_3)\otimes V(2_63_3)\otimes V(2_0)\otimes V(2_6)\to V(\bs\pi)\otimes V(2_6)\\ V(2_6^23_3)\otimes V(1_32_0)\otimes V(1_3)\to V(\bs\pi)\otimes V(1_3).
	\end{gather*}
	In particular, all the ordered two-fold tensor products in 
	\begin{equation*}
		V(\bs\pi)\otimes V(2_6)\otimes V(1_3) 
	\end{equation*}
	are highest-$\ell$-weight, so the whole tensor product is highest-$\ell$-weight by \Cref{t:cyc}.  By considering the epimorphism
	\begin{equation*}
		V(\bs\pi)\otimes V(2_6)\otimes V(1_3)\to V(\bs\pi)\otimes V(\bs\pi_H),
	\end{equation*}
	we complete the proof of \eqref{e:realbutnotstrongiii}. 
	
	For illustrative purposes, let us give an alternative proof that $H$ satisfies condition (iii) of \Cref{d:rds}. Begin by noting that 
	\begin{equation}\label{e:13263326flip}
		V(2_63_3)\otimes V(\bs\pi_H) \quad\text{is highest-$\ell$-weight.}
	\end{equation}
	This can be proved with an argument used to prove \eqref{e:13263326} (or use the symmetry of $G(1_33_32_6^2)$). Hence, it follows from \Cref{p:vnvstar}, \eqref{e:13263326}, \eqref{e:13263326flip}, and \Cref{kkop}(a) that 
	\begin{equation*}
		\lie d (V(\bs\pi_H), V(3_32_6)) = 0.
	\end{equation*}
	Then, by \eqref{e:kkopineq}  and Lemma \ref{l:naoi}, we have
	\begin{equation*}
		\lie d (V(\bs\pi_H), V(\bs\pi_{H^c}))\leq \lie d (V(\bs\pi_H), V(3_32_6))+ \lie d (V(\bs\pi_H), V(2_0)) = \lie d (V(\bs\pi_H), V(2_0))=1.
	\end{equation*}
	Therefore, \Cref{c:kkop} implies 
	\begin{equation*}
		V(\bs\pi)\otimes V(\bs\pi_H), \ \ \ V(\bs\pi)\otimes V(\bs\pi_{H^c}) \quad\text{are simple},
	\end{equation*}
	which proves that both $H$ and $H^c$ satisfy condition (iii) of \Cref{d:rds}. Since conditions (i) and (ii) are symmetric on $H$ and $H^c$, it follows that both $H$ and $H^c$ are rds in this case. However, it is not always true that $H^c$ is an rds if $H$ is, as we shall see in \Cref{ex:realbutnotstrongrev}.
	
	Since $H^c$ is a tree, it is strongly real by \Cref{t:srtree} below. In fact, the theorem implies that $2_6, 3_3, 2_0$ (or the opposite order) is an rds-quochain for $H^c$ and, hence, $H, 2_6, 3_3, 2_0$ is an rds-quochain for $G$ of length $4$. Similarly, $H^c, 2_6,1_3$ is  a length-3 rds-quochain for $G$. We will next show that $V(\bs\pi)$ is not strongly real, which means there is no rds-quochain for $G$ of lengths $5$. Since $G$ is real, it is, by itself, a weak rds-quochain of length $1$.  
					
	We will show that $V(\bs\pi)$ is not strongly real by showing that no vertex of $G(\bs\pi)$ is an rds. Begin by noting that the sources of $G(\bs\pi)$ do not determine real cuts, since $V(\bs\pi2_6^{-1})$ is imaginary by \cite{lec:im}. For the vertex $1_3$, an application of \eqref{e:sJs} with $J=[1,2]$ implies
	\begin{equation*}
		V(1_3)\otimes V(\bs\pi1_3^{-1})\quad\textrm{and}\quad V(\bs\pi1_3^{-1})\otimes V(1_3) \ \ \text{are not highest-$\ell$-weight},
	\end{equation*}
	so $1_3$ does not satisfy condition (ii) for being an rds. Indeed,
	\begin{equation*}
		V(1_3)_J\otimes V(\bs\pi1_3^{-1})_J\cong 	V(1_3)\otimes_J V(2_02_6^2)\cong V(1_3)\otimes_J V(2_6)\otimes_J V(2_02_6) 
	\end{equation*}
	and $V(1_3)\otimes_J V(2_6)$ is not highest-$\ell$-weight. On the other hand,
	\begin{equation*}
		V(\bs\pi1_3^{-1})_J\otimes V(1_3)_J \cong 	V(2_02_6^2)\otimes_J V(1_3)\cong V(2_6^2)\otimes_J V(2_0)\otimes_J V(1_3)
	\end{equation*}
	and $V(2_0)\otimes_J V(1_3)$ is not highest-$\ell$-weight.
	Similarly, $3_3$ is also not an rds.
				
	It remains check $2_0$ is not an rds as well.  We will show $V(2_0)\otimes V(\bs\pi)$ is not highest-$\ell$-weight, so condition (iii) of \Cref{d:rds} is not satisfied. Indeed, if this were the case, since $V(2_6^2)\otimes V(\bs\pi2_6^{-2})$ is highest-$\ell$-weight, \Cref{p:hlwmorph}(ii) would imply that 
	\begin{equation*}
		V(2_0)\otimes V(\bs\pi2_6^{-2})=V(2_0)\otimes V(1_33_32_0)\quad\text{is highest-$\ell$-weight}.
	\end{equation*}
	However, since $V(1_3)\otimes V(3_32_0)$ is highest-$\ell$-weight and ${~^*}V(2_0)\otimes V(3_32_0)=V(2_4)\otimes V(3_32_0)$ is simple by \Cref{t:3lineprime}, \Cref{p:hlwdualmorph} would imply that $V(2_0)\otimes V(1_3)$ is highest-$\ell$-weight, yielding a contradiction. For completeness, let us check the assumptions of \Cref{t:3lineprime} are indeed satisfied. In the notation of that theorem, let $\bs\omega_{i_j} = 2_4, \bs\omega_{i_{j'}}=3_3$, so $m_j=4,m_{j'}=3$, and $I_j=I$. The last of these facts immediately implies the first 2 conditions of the theorem are satisfied. The third condition becomes $3-4+4\in\mathscr R_{2,3}$, which is true, while \eqref{2exracondm} reads $4+1\le 3+1+1$, which is also true. \endd
\end{ex}

Let us streamline one argument used  in \Cref{ex:realbutnotstrong}. To shorten notation, given graphs $H$ and $K$, we set
\begin{equation}\label{e:kkopgraph}
	\lie d(H,K) = \lie d(V(\bs\pi_H),V(\bs\pi_K)). 
\end{equation} 

\begin{lem}\label{l:rdsbysmallkkop}
	Suppose $H$ satisfies conditions (i) and (ii) of \Cref{d:rds}. If $\lie d(H,H^c)\le 1$, then both $H$ and $H^c$ are rds for $G$.
\end{lem}

\begin{proof}
	Since conditions (i) and (ii) are symmetric on $H$ and $H^c$, it remains to show 
	\begin{equation*}
		\lie d(H,H^c)\le 1 \quad\Rightarrow\quad V(\bs\pi_G)\otimes V(\bs\pi_H) \ \ \text{and}\ \ V(\bs\pi)\otimes V(\bs\pi_{H^c}) \ \ \text{are simple}. 
	\end{equation*}
	By renaming if necessary, we can assume condition (ii) is
	\begin{equation*}
		V(\bs\pi_H)\otimes V(\bs\pi_{H^c}) \quad\text{is highest-$\ell$-weight}.
	\end{equation*}
	In particular, $V(\bs\pi_G)={\rm hd}(V(\bs\pi_H)\otimes V(\bs\pi_{H^c}))$. If $\lie d(H,H^c)= 1$, \Cref{c:kkop} completes the proof since both $V(\bs\pi_H)$ and $V(\bs\pi_{H^c})$ are real by condition (i). If $\lie d(H,H^c)= 0$, then $V(\bs\pi_G)\cong V(\bs\pi_H)\otimes V(\bs\pi_{H^c})\cong V(\bs\pi_{H^c})\otimes V(\bs\pi_H)$ and we get
	\begin{gather*}
		V(\bs\pi_G)\otimes V(\bs\pi_H)\cong V(\bs\pi_{H^c})\otimes V(\bs\pi_H)\otimes V(\bs\pi_H),\\
		V(\bs\pi_G)\otimes V(\bs\pi_{H^c})\cong V(\bs\pi_{H})\otimes V(\bs\pi_{H^c})\otimes V(\bs\pi_{H^c}).
	\end{gather*}
	The assumptions imply all ordered two-fold tensor products are simple, and we are done by \Cref{c:hlwquot}.
\end{proof}

\begin{ex}\label{ex:realbutnotstrongrev}
    Let us return to \Cref{ex:realbutnotstrong}. We have shown there that $2_0$ does not satisfy condition (iii) of \Cref{d:rds}. Note that it satisfies conditions (i) and (ii). Indeed, since $2_0$ is a sink, $V(\bs\pi2_0^{-1})\otimes V(2_0)$ is highest-$\ell$-weight, showing condition (ii) is satisfied. For condition (i), we have checked in  \Cref{ex:realbutnotstrong} that $V(2_61_3)\otimes V(2_63_3)$ is simple and, hence, $V(\bs\pi2_0^{-1})\simeq V(2_61_3)\otimes V(2_63_3)$ is real since $G(1_32_6)$ and $G(3_32_6)$ are trees and, therefore, real by \cite[Theorem 2.4.8]{ms:3tree}. 
	
	It then follows from \Cref{l:rdsbysmallkkop} that $\lie d(2_0,G(1_33_32_6^2))>1$. It can be shown that it is actually $2$. The role of pairs of modules whose KKOP invariant is larger than $1$ has not been systematically studied so far. We want to give an example of an rds $H$ such that $H^c$ is not an rds and we believe $G(1_33_32_6^2)$ is not an rds for $G$, which would provide such an example. However, the computations for checking this with the techniques we have employed are inconclusive for deciding if that is really the case. On the other hand, letting $\bs\varpi=1_31_5^22_02_22_8$ for type $A_n, n\ge 2$,  then $G=G(\bs\varpi)$ is
	
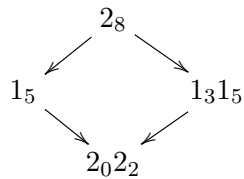
\begin{equation*}
		\xymatrix@R-1pc@C-1pc{ & 2_8 \ar[ld] \ar[rd] &  \\	1_5 \ar[rd] & & 1_31_5 \ar[ld] \\ & 2_02_2 & }
	\end{equation*}
	In this case, we have checked that $H=G(\bs\varpi 2_8^{-1})$ is an rds, but $H^c$ is not. Since the computations are rather lengthy, we have chosen not to include them here.	\endd
\end{ex}

\subsection{Valence-1 Vertices}\label{ss:eas}
The purpose of this section is to describe a sufficient condition for a valence-1 vertex of a pseudo $q$-factorization graph to be an rds. We start by streamlining one of the key steps of the proof of \cite[Theorem 2.4.8]{ms:3tree}, which was originally proved by a systematic use of \Cref{t:cyc} and its corollaries, alongside \Cref{p:vnvstar}. Here, we give an alternative proof using KKOP invariants, which lead to a shorter proof for a stronger statement. 

\begin{prop}\label{p:eas}
	Let $\bs\omega$ be an extremal vertex of a pseudo $q$-factorization graph $G$ over $\bs\pi$ and $\bs\varpi = \bs\pi_{\Adj(\bs\omega)}$. If $V(\bs\varpi)\otimes V(\bs\omega)$ is simple, then $V(\bs\pi)\otimes V(\bs\omega)$ is simple. In particular, if $\bs\omega$ determines a real cut of $G$, then $\bs\omega$ is an rds of $G$.
\end{prop}

\begin{proof}
	We freely use the fact that $V(\bs\omega)$ is real when using \Cref{kkop} in what follows. Note $G\setminus \Adj(\bs\omega)$ and $\bs\omega$ lie in distinct connected components of $(G\setminus \Adj(\bs\omega))\otimes \{\bs\omega\}$. In particular, \Cref{t:cyc} implies
	 \begin{equation}\label{adjout}
	 	V(\bs\omega)\otimes V(\bs\pi\bs\varpi^{-1}) \ \ {\rm is\ simple}.
 	\end{equation} 
	By \eqref{e:kkopineq} we have
	\begin{equation*}
		\lie d (V(\bs\omega), V(\bs\pi))\leq \lie d (V(\bs\omega),V(\bs\varpi)) + \lie d (V(\bs\omega),V(\bs\pi\bs\varpi^{-1})).
	\end{equation*}
	The assumption that  $V(\bs\varpi)\otimes V(\bs\omega)$ is simple and \eqref{adjout}, together with \Cref{kkop}(a), imply each of terms in the right hand side is zero and so $\lie d(V(\bs\omega), V(\bs\pi)) = 0$. Another application of Proposition \ref{kkop}(a) completes the proof of the first statement. To prove the "in particular", it suffices to show that condition (ii) of \Cref{d:rds} holds. But this follows from \Cref{l:topbothw}, since $\bs\omega$ is extremal.    
\end{proof}

The ``in particular'' part of the next corollary recovers \cite[Theorem 2.4.8]{ms:3tree}.

\begin{cor}\label{c:eas2}
	Suppose $\lie g$ is of type $A$ and that $\bs\omega$ determines a real cut of a pseudo $q$-factorization graph $G$. If  $\bs\omega$ has valence at most $1$ in  $G$, then $\bs\omega$ is an rds in $G$. In particular, if $G$ is a tree, $G$ is strongly real.
\end{cor}
			
\begin{proof}
	If $|\hat A_{\bs\omega}|=0$, then $V(\bs\pi_G)\cong V(\bs\omega)\otimes V(\bs\pi\bs\omega^{-1})$ by \eqref{e:dissosimple}, which implies \eqref{e:realts}. Thus, assume  $|\hat A_{\bs\omega}|= 1$. In this case, $\bs\pi_{\Adj(\bs\omega)} = \bs\omega\bs\varpi$ for some $\bs\varpi\in\mathcal{KR}$. Then, by \Cref{p:eas}, in order to show $\bs\omega$ is an rds, we are left to show $V(\bs\omega\bs\varpi)\otimes V(\bs\omega)$ is simple. But this is \Cref{l:3aline}. The in particular part follows by an easy induction on the number of vertices of $G$. Indeed, if $G$ is a tree, then $G$ has a valence-1 vertex, say $v$, and $G\setminus v$ is also a tree. Hence, $v$ is extremal and determines a real cut (by inductive assumption), so conditions (i) and (ii) of \Cref{d:rds} are satisfied and condition (iii) follows by  the first statement.
\end{proof}
			
\begin{rem}\label{r:treesr}
	The reason we had to assume that $\lie g$ is of type $A$ in \Cref{c:eas2} is because we used \Cref{l:3aline}.  However, although we are not aware of a published proof of this fact for more general $\lie g$, it is believed that its statement is true as long as the minimal subdiagram of $I$ on which all the structure of corresponding $q$-factorization graphs make sense does not contain a trivalent node. In that case, the above proof is type independent. \Cref{l:3aline} is known to be false in the presence of the trivalent node.\endd
\end{rem}

\begin{ex}
	Let us illustrate how the above test depends on the chosen pseudo $q$-factorization. Let $\lie g$ be of type $A$ with $n\ge 7$ and $\bs\pi = 1_5 2_2 4_04_24_4 4_8 6_2 7_5$. Then, the following are pseudo $q$-factorization graphs over $\bs\pi$.
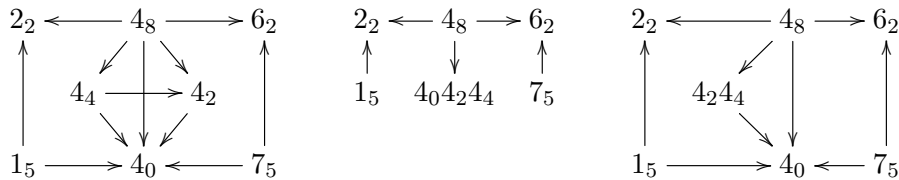
	\begin{equation*}
		\xymatrix@R-1pc@C-1.5pc{2_2  &  & 4_8\ar[ll]\ar[rr]\ar[dd]\ar[dl]\ar[dr] & & 6_2 \\ & 4_4\ar[rr]\ar[dr] & & 4_2\ar[dl] & \\ 1_5\ar[uu]\ar[rr] & & 4_0 & & 7_5\ar[ll] \ar[uu]} \qquad 
		\xymatrix@R-1pc@C-1.5pc{2_2 & 4_8\ar[d]\ar[l]\ar[r]  & 6_2 \\ 1_5\ar[u] & 4_04_24_4 & 7_5\ar[u]}
		\qquad \xymatrix@R-1pc@C-1.5pc{2_2 & & 4_8\ar[dd]\ar[dl]\ar[rr]\ar[ll] & & 6_2\\& 4_24_4\ar[dr]  \\ 1_5\ar[rr]\ar[uu] & & 4_0 & &7_5\ar[ll]\ar[uu]}
	\end{equation*}
	The first is the fundamental factorization graph and the middle is the actual $q$-factorization graph. The adjacency graphs of all the sources in all three graphs are extremal. However, \Cref{c:eas2} is not applicable for the first and last graphs because $|\hat A_{\bs\omega}|\ge 2$ for all sources. Since the middle graph is a tree, \Cref{c:eas2} is applicable and, hence, $V(\bs\pi)$ is real. If we insisted on using the first or last graphs and we chose $\bs\omega=1_5$, then the adjacency subgraph is $\xymatrix@R-1pc@C-1pc{2_2 & 1_5\ar[r]\ar[l] & 4_0}$ and we would need to check whether $V(1_52_24_0)\otimes V(1_5)$ is simple. 
	Using cluster algebra arguments, it is possible to show  $V(1_52_24_0)\cong V(1_52_2) \otimes V(4_0)$ (see, for instance, \cite[Proposition 3.5]{bc:tphl}), which implies $V(1_52_2)\otimes V(4_0)\otimes V(1_5)$ is reducible, thus showing that $\bs\omega$ is not an rds.\endd 
\end{ex}
			
\subsection{Prime Snakes as Adjacency Subgraphs}\label{ss:etoas}
We now prove that an extremal vertex whose adjacency subgraph corresponds to a prime snake is an rds.  We assume throughout this section that  $\lie g$ is of type $A$.
The class of prime snake modules was originally introduced in \cite{muyou:tsystem} (see also \cite{DLL:snakes}). Having the concept of fundamental factorization graph in mind, one immediate notices that the definition is equivalent to requiring that the corresponding fundamental factorization graph is totally ordered.   For instance, in \Cref{l:naoi}, the assumption on $\bs\pi$ can be rephrased as ``$V(\bs\pi)$ is a prime snake module''.   The proof of the following alternate characterization will appear in \cite{bms:tof}. 

\begin{lem}\label{l:psiffcon}
	Let $\bs\pi\in \mathcal P^+$ be such that $V(\bs\pi)$ is a snake module and  let $G$ be a pseudo $q$-factorization graph of $\bs\pi$. The following are equivalent: 
	\begin{enumerate}[(i)]
		\item $V(\bs\pi)$ is prime;
		\item $G$ is totally ordered;
		\item $G$ is connected.\qed
	\end{enumerate}
\end{lem}

The most crucial previously proved result for this subsection is the following.
	
\begin{prop}[{\cite[Proposition 4.1.3(ii)]{naoi:Tsys}}]\label{p:testworksforsnakes}
	Suppose $V(\bs\pi)$ is a prime snake module and that $\bs\omega$ divides $\bs\pi$. Then, $V(\bs\pi)\otimes V(\bs\omega)$ is simple.  \qed
\end{prop}
			
Note \Cref{p:testworksforsnakes} implies, in particular, that prime snake modules are real. 			
			
\begin{cor}\label{c:testworksforsnakes}
	Suppose $\bs\omega$ is an extremal vertex of a pseudo $q$-factorization graph $G$ which determines a real cut.	If $V(\bs\pi_{\Adj(\bs\omega)})$ is a prime snake module, $\bs\omega$ is an rds in $G$. 
\end{cor}
			
\begin{proof}
	\Cref{p:testworksforsnakes} implies all the assumptions of  \Cref{p:eas} hold. 
\end{proof}

\begin{cor}\label{c:snakesarereal}
	Suppose $V(\bs\pi)$ is a prime snake module and that $\bs\varpi$ is any vertex in $G$. Then, there exists an rds-quochain, say $\mathcal G=G_1,\dots,G_l$, such that $G_k$ is an extremal vertex of $\bar G_{k-1}$ for all $1\le k\le l$ and $\bs\varpi$ is the vertex of $G_l$. In particular, $G$ is strongly real.
\end{cor}
			
\begin{proof}
	We proceed by induction on the number of vertices $l$ of $G$, which clearly starts when $l=1$. Thus, assume $l>1$ and let $\bs\omega_1,\dots,\bs\omega_l$ be an enumeration of the vertices of $G$ compatible with the total order of $G$ and so that $\bs\omega_1$ is the source. If $\bs\varpi=\bs\omega_1$, let $\bs\omega=\bs\omega_l$ and, otherwise, let $\bs\omega=\bs\omega_1$. We will show that the subgraph $G_1$ having $\bs\omega$ as its unique vertex is an rds. Since $G\setminus G_1$ is a pseudo $q$-factorization graph over $\bs\pi\bs\omega^{-1}$ and $V(\bs\pi\bs\omega^{-1})$ is a prime snake module, the induction hypothesis implies there exists an rds-quochain $G_2,\dots, G_l$  for $G\setminus G_1$ satisfying the properties in the statement. It follows that $G_1,\dots, G_l$ is an rds-quochain for $G$ as desired. 
				
	Since $\Adj(\bs\omega)$ is connected by definition, it follows from \Cref{l:psiffcon} that $V(\bs\pi_{\Adj(\bs\omega)})$ is a prime snake module.
	Moreover, $\Adj({\bs\omega})$ is an extremal subgraph of $G$. The induction hypothesis implies, in particular, that $V(\bs\pi\bs\omega^{-1})$ is real, so $\bs\omega$ determines a real cut. It then follows from \Cref{c:testworksforsnakes} that $G_1$ is an rds as claimed, thus completing the proof.
\end{proof}

\section{Applications}\label{s:applications}

\subsection{Multicuts of Tree Type}\label{ss:srt}
We now describe a procedure for constructing new strongly real modules from certain known ones. We start with a generalization of the notion of a tree.
			
Recall the definition of tensor products of pseudo $q$-factorization graphs given in the paragraph containing \eqref{e:dissosimple}. Let $\mathcal G$ be a finite sequence of length $l$ of pseudo $q$-factorization graphs, say $G_k=(\mathcal V_k,\mathcal A_k), 1\le k\le l$ and set $G:=G_1\otimes G_2\otimes\cdots\otimes G_l$. In particular, $\mathcal G$ is a multicut of $G$. Recall \eqref{e:cutsetgen}. We say $G$ is a $\mathcal G$-tree if $G$ is connected and 
\begin{equation}
		\#\mathcal A_{\mathcal G} = l-1.
\end{equation}
Note that, if $G_k$ has a single vertex and no loops for all $k$, this is one of the usual definitions of a tree.
Given $k\ne k'$, let us say $G_k$ is linked to $G_{k'}$ if there exists $a\in\mathcal A_{G}$ such that $h(a)\in G_k$ and $t(a)\in G_{k'}$ or vice-versa. Then,   define the valence of $G_k$ in $G$ as
\begin{equation*}
	{\rm val}_G(G_k) = \#\{k'\ne k: G_k \text{ is linked to } G_{k'}\}.
\end{equation*}
As in the case of usual trees, if $\ell>1$,  there exists at least two values of $k$ such that ${\rm val}(G_k)=1$. Indeed, we can form a graph $TG$ whose vertex set is $\{G_k:1\le k\le l\}$ and the arrow set $\mathcal A_{TG}$ is determined by
\begin{equation*}
	(G_k,G_m)\in\mathcal A_{TG} \quad\Leftrightarrow\quad \exists\ v\in\mathcal V_{G_k}, \ w\in\mathcal V_{G_m} \ \ \text{such that}\ \ (v,w)\in \mathcal A_{G}.
\end{equation*}
One easily checks there is a bijection $\mathcal A_{\mathcal G}\to\mathcal A_{TG}$ and
\begin{equation*}
	{\rm val}_G(G_k) = {\rm val}_{TG}(G_k),
\end{equation*}
where the latter is the usual valence of a vertex of a graph. In particular, we can always assume, up to ordering the sequence $\mathcal G$, that 
\begin{equation}\label{e:val1reo}
	{\rm val}_{\bar G_{k-1}}(G_k) =1 \quad\text{for all}\quad 1\le k<l. 
\end{equation}
			
We will say $G$ can be realized as a snake tree if there exists a multicut $\mathcal G$ of $G$ such that $G$ is a $\mathcal G$-tree and $V(\bs\pi_{G_k})$ is a snake module for all $k$.  Note that if $G$ is a tree, then it can be realized as a snake tree (with $\mathcal G$ being any enumeration of the vertices). More generally, if $G_k$ satisfies a certain property $P$, we shall say $G$ is a $P$ tree (with parts in $\mathcal G$).  
			
\subsection{Snake Trees are Strongly Real}\label{ss:snaketrees}
Let us start with a couple of examples which illustrate the general argument for showing that snake trees are strongly real, a fact that follows from the more general statement of \Cref{t:srtree} below.
			
\begin{ex}\label{ex:prst}
	Let $\lie g$ be of type $A$ with $n\ge 6$ and $\bs\pi = \bs\omega_{4,0}\bs\omega_{2,4}\bs\omega_{3,9,3}\bs\omega_{2,14,3}$. Note this is the $q$-factorization and associated $q$-factorization graph is 
	\begin{equation*}
		\xymatrix@R-1pc@C-1pc{ \bs\omega_{4,0} & \ar[l]\bs\omega_{2,4} &\\ & \ar[u]\ar[ul]\bs\omega_{3,9,3} & \ar[l]\bs\omega_{2,14,3} }
	\end{equation*}
	Indeed,
	\begin{gather*}
		\mathscr R_{4,2} = \{4,6\}, \quad \mathscr R_{4,3}^{1,3} = \{5,7,9\}, \quad \mathscr R_{4,2}^{1,3} = \{6,8\}, \\
		\mathscr R_{2,3}^{1,3} = \{5,7\}, \quad \mathscr R_{2,2}^{1,3} = \{4,6\}, \quad \mathscr R_{2,3}^{3,3} = \{3,5,7,9\}.
	\end{gather*}
	Then, letting $G_1 = G(\bs\omega_{4,0}\bs\omega_{2,4}\bs\omega_{3,9,3})$ and $G_2=G(\bs\omega_{2,14,3})$, it follows that $G(\bs\pi)=G_1\otimes G_2$ is a snake tree, since $G_1$ and $G_2$ are snake modules. One easily checks $V(\bs\pi)$ is not a snake module nor there exists a pseudo $q$-factorization over $\bs\pi$ which is a tree. Note that \Cref{c:eas2} applies to the unique source, which is the unique vertex of $G_2$,  so $G_2$ is an rds. After ``quotienting out'' $G_2$, we remain with $G_1$, which is a prime snake, so the original module is strongly real.  		
  
   Let us make a remark aimed at comparing the results of this section with those of  \Cref{ss:kkkords} below. A repeated use of \Cref{kkop} and \Cref{l:naoi} allows us to deduce that 
 	$$\lie d (V(\bs\pi_{G_1}), V(\bs\pi_{G_2}))\leq 2.$$
 	Since $G$ is totally ordered, it is prime by \cite[Theorem 3.5.5]{ms:to} and, hence, $\lie d (V(\bs\pi_{G_1}),V(\bs\pi_{G_2}))>0$ by \Cref{kkop}(a). We believe $\lie d (V(\bs\pi_{G_1}),V(\bs\pi_{G_2}))=2$, but we could not find a way of showing it is at least $2$ using known results, except, possibly, by studying the order of the zero of the polynomial $d_{G_1,G_2}(z)d_{G_2,G_1}(z)$ at $z=1$, where $d_{G_i,G_j}(z)$ denotes the denominator of the normalized $R$-matrix corresponding to $V(\bs\pi_{G_i})$ and $V(\bs\pi_{G_j})$. In any case, this shows that $\lie d (V(\bs\pi_{G_1}),V(\bs\pi_{G_2}))$ plays no role in the above argument for showing $G$ is strongly real, in contrast to \Cref{t:kkoprds}   (see the paragraph of the introduction dedicated to explaining the differences between \Cref{t:kkoprds} and \cref{t:srtree}). 	\endd
\end{ex}

\begin{ex}\label{ex:nextglue}
	Let $\lie g$ be of type $A$ with $n\ge 6$ and $\bs\pi = \bs\omega_{4,0}\bs\omega_{2,4}\bs\omega_{3,9,3}\bs\omega_{1,7}$. Note this is the $q$-factorization and the associated $q$-factorization graph is 
	$\xymatrix@R-1pc@C-1pc{ \bs\omega_{4,0} & \ar[l]\bs\omega_{2,4} &  \ar[l]\bs\omega_{1,7} \\ & \ar[u]\ar[ul]\bs\omega_{3,9,3} }$.
	This is a snake tree with $G_1 = G(\bs\omega_{4,0}\bs\omega_{2,4}\bs\omega_{3,9,3})$ and $G_2=G(\bs\omega_{1,7})$, but it is not a snake. The results from Sections \ref{ss:eas} and \ref{ss:etoas} 
	imply all extremal vertices are rds for $G$. Indeed, \Cref{c:eas2} implies $\bs\omega_{1,7}$ is an rds, while \Cref{c:testworksforsnakes} implies the same  holds for the other two. In the case of $\bs\omega_{1,7}$, what remains is a snake module, while for the others what remains is a tree. In any case, we conclude the original module is strongly real. 	\endd
\end{ex}

We can actually consider a larger class of $\mathcal G$-trees, which we now define.

\begin{defn}
	We say that a vertex $v$ of a pseudo $q$-factorization graph $G$ is an rds-base for $G$ if there exists an rds-quochain having $v$ as the final rds. If there is such a quochain so that all the rds are single vertices, then $v$ is said to be a strong rds-base. 	
	A $\mathcal G$-tree is said to be well-based if, for every $a\in\mathcal A_{\mathcal G}$, both $t(a)$ and $h(a)$ are rds-bases for the corresponding subgraphs of the multicut $\mathcal G$. If both $t(a)$ and $h(a)$ are strong rds-bases for every $a\in\mathcal A_{\mathcal G}$, the $\mathcal G$-tree is said to be strongly based.  \endd
\end{defn}

\Cref{c:snakesarereal} implies all pseudo $q$-factorization graphs associated to prime snake modules are fully strongly based, i.e., every vertex is a strong rds-base. Evidently, if $\mathcal G=G_1,\dots, G_l$ with $G_k$ fully (strongly) based for all $k$, then $\mathcal G$ is well (strongly) based. In particular, if $G_k$ is a singleton for every $k$, i.e., if $G$ is a usual tree, then $\mathcal G$ is strongly based. 
			
Henceforth, assume \eqref{e:val1reo} holds and that $\mathcal G_k = G_{k,1},\dots, G_{k,m_k}$ is an rds-quochain for $G_k$. 
In particular, 
\begin{equation*}
	G_{k,m_k} = \{v_k\} \ \text{for some}\ \ v_k\in\mathcal V_{G_k} \quad\text{and}\quad
	G_k = G_{k,1}\otimes\cdots\otimes G_{k,m_k} \quad\text{for all}\quad 1\le k\le l.
\end{equation*}
We shall denote by $\mathcal G_1*\cdots*\mathcal G_l$ the resulting concatenation of the rds-quochains $\mathcal G_1,\dots,\mathcal G_l$. In general, there is no reason for it to be an rds-quochain for $G$. We shall say $\mathcal G_k$ is compatible with $\mathcal G$ if
\begin{equation}\label{e:compat}
	v_k \in\{t(a):a\in\mathcal A_{\mathcal G}\}\cup \{h(a):a\in\mathcal A_{\mathcal G}\}.
\end{equation}
Note that, if $G$ is a tree, i.e. if $G_k$ is a singleton for all $k$, then $m_k=1$ for all $k$ and $\mathcal G_k$ is compatible with $\mathcal G$. Thus, if $G$ is a tree, all the assumptions of the following theorem are vacuously satisfied. On the other hand, since any graph associated to a snake module is fully strongly based, there always exists a collection of rds-quochains $\mathcal G_k$ which are  compatible with $\mathcal G$.  

\begin{ex} Consider $G= G_1\otimes G_2$ as in \Cref{ex:nextglue}, set $\mathcal G=G_1, G_2$ and note that $\mathcal A_{\mathcal G} =\{a\}$ where $a$ is 
\begin{equation*}
	\xymatrix{ \bs\omega_{2,4}  & \ar[l]_{a} \bs\omega_{1,7}}.
\end{equation*}	
In particular, \eqref{e:val1reo} is trivially satisfied. Since $G_2$ has $\bs\omega_{1,7}$ as its unique vertex, it is evidently an rds-base for $G_2$. On the other hand, $G_1$ is a snake module, so every vertex is a strong rds-base for $G_1$ and, hence, $\mathcal G$ is strongly based. 	
	
Setting $G_{1,1} = G(\bs\omega_{4,0}), G_{1,2} = G(\bs\omega_{3,9,3}), G_{1,3}= G(\bs\omega_{2,4})$, it follows that $\mathcal G_1 = G_{1,1}, G_{1,2}, G_{1,3}$ is an rds-quochain for $G_1$ ending at the vertex $\bs\omega_{2,4}$. Since $G_2$ is a singleton, it is, by itself, its unique rds-quochain. To match the general notation above, set $\mathcal G_2 = G_{2,1}$ with $G_{2,1}=G_2$. Then, $v_1 = \bs\omega_{2,4}, v_2 = \bs\omega_{1,7}$, and it follows that \eqref{e:compat} is satisfied, i.e., $\mathcal G_1,\mathcal G_2$ are compatible with $\mathcal G$. Moreover, $\mathcal G_1*\mathcal G_2 = G_{1,1}, G_{1,2}, G_{1,3},  G_{2,1}$ and the comments made in \Cref{ex:nextglue} show it is an rds-quochain for $G$. \endd
\end{ex}
            
\begin{thm}\label{t:srtree}
	Let $\lie g$ be of type $A$. Suppose $\mathcal G$ is well-based, satisfies \eqref{e:val1reo}, and $\mathcal G_k$ is an rds-quochain compatible with $\mathcal G$ for all $1\le k\le l$. Then,  $\mathcal G_1*\cdots*\mathcal G_l$  is an rds-quochain for $G$. In particular, if $m_k = \#\mathcal V_{G_k}$ for all $1\le k\le l$, then $G$ is strongly real. 
\end{thm}
			
\begin{proof}
	Let us proceed by induction on the number $N$ of vertices of $G$. Since there is nothing to do if $l=1$, induction starts. Thus, assume $l>1$ and let 
	\begin{equation*}
		G_1' = G_1\setminus G_{1,1}, \quad \mathcal G' = G_1',G_2,\dots,G_l, \quad\text{and}\quad G'  = G\setminus G_{1,1} = G_1'\otimes G_2\otimes\cdots\otimes G_l.		
	\end{equation*} 
	In order to complete the proof we need to check
	\begin{equation}\label{e:srtree}
		G_{1,1} \ \ \text{is an rds for $G$ and}\ \ \mathcal G' \ \ \text{satisfies the inductive assumption.}
	\end{equation}
	
	Let us start with the second part of \eqref{e:srtree}. 	
	If $G_1=G_{1,1}$, then $G_1'$ is the empty graph (so it is omitted above) and $\mathcal G'$ obviously satisfies all the assumptions imposed on $\mathcal G$. Otherwise, $v_1$ is a vertex of $G_1'$ and it is the unique vertex linking it to another part of $\mathcal G'$. This immediately implies $\mathcal G'$ is well-based and satisfies \eqref{e:val1reo}. Moreover, the sequence $\mathcal G_1' = G_{1,2},\dots,G_{1,m_1}$ is an rds-quochain for $G_1'$ ending in $v_1$, so it is compatible with $\mathcal G'$, thus completing the checking of the second part of \eqref{e:srtree}. Let us turn to the first.
	
	Suppose first that $G_1 = G_{1,1}$ and, hence, $\mathcal V_{G_1}=\{v_1\}$ since $\mathcal G_1$ is an rds-quochain.  By inductive assumption, $G'$ is real, showing that $G_{1,1}$ satisfies part (i) of \Cref{d:rds}. Since $G_1$ has valence $1$, \Cref{c:eas2} implies it is an rds for $G$. 	
	
	Henceforth, assume $m_{1,1}>1$, so $v_1$ is not a vertex of $G_{1,1}$. Since $G_{1,1}$ is an rds for $G_1$, it is real and the inductive assumption implies $G\setminus G_{1,1}=G'$ is real, so $G_{1,1}$ determines a real cut in $G$. For checking condition (ii) of \Cref{d:rds}, using arrow duality if necessary, assume 
	\begin{equation*}
		\Adj_G(v_1) = v_1\longrightarrow v \ \ \text{for some vertex $v$ of}\ \ G_2
	\end{equation*} 
	 and set
	\begin{equation*}
		G'' = G_2\otimes \cdots\otimes  G_l.
	\end{equation*}
	Since $G_1$ has valence $1$, together with \Cref{c:hlwpqf}, this implies
	\begin{equation}\label{e:srtreenhwH}
		V(\bs\pi_H)\otimes V(\bs\pi_{G''}) \quad\text{is highest-$\ell$-weight for all}\quad H\triangleleft G_1,
	\end{equation}
	and, moreover,
	\begin{equation}\label{e:srtreenov1}
		v_1\notin\mathcal V_H \quad\Rightarrow\quad V(\bs\pi_H)\otimes V(\bs\pi_{G''}) \ \ \text{is simple}.
	\end{equation}
	Using that $G_{1,1}$ satisfies \Cref{d:rds}(ii) when regarded as a subgraph of $G_1$, there are two possibilities:
	\begin{equation}\label{e:srtree2p}
		V(\bs\pi_{G_{1,1}})\otimes V(\bs\pi_{G_1'}) \ \ \text{is highest-$\ell$-weight or}\ \  	V(\bs\pi_{G_1'})\otimes V(\bs\pi_{G_{1,1}})\ \ \text{is highest-$\ell$-weight.}
	\end{equation}
	In the first case, all the ordered two-fold tensor products arising from the following tensor product are highest-$\ell$-weight
	\begin{equation*}
			V(\bs\pi_{G_{1,1}})\otimes V(\bs\pi_{G_1'})\otimes V(\bs\pi_{G''}). 
	\end{equation*}
	Moreover, the three factors are real modules. Hence, \Cref{t:cyc} implies this tensor product is highest-$\ell$-weight.  Since and $G'=G_1'\otimes G''$, we also have an epimorphism
	\begin{equation*}
			V(\bs\pi_{G_{1,1}})\otimes V(\bs\pi_{G_1'})\otimes V(\bs\pi_{G''}) \to 	V(\bs\pi_{G_{1,1}})\otimes V(\bs\pi_{G'}),
	\end{equation*}
	which shows \Cref{d:rds}(ii) holds in this case since it follows that
	\begin{equation*}
			V(\bs\pi_{G_{1,1}})\otimes V(\bs\pi_{G'}) \quad\text{is highest-$\ell$-weight.}
	\end{equation*} 
	If the second possibility in \eqref{e:srtree2p} holds, consider instead 
	\begin{equation*}
		V(\bs\pi_{G_1'})\otimes V(\bs\pi_{G''})\otimes V(\bs\pi_{G_{1,1}}). 
	\end{equation*}
	This time,  we use \eqref{e:srtreenov1} to conclude $V(\bs\pi_{G''})\otimes V(\bs\pi_{G_{1,1}})$ is highest-$\ell$-weight and similar arguments then imply 
	\begin{equation*}
		V(\bs\pi_{G'}) \otimes  V(\bs\pi_{G_{1,1}}) \quad\text{is highest-$\ell$-weight,}
	\end{equation*} 
	thus completing the poof that \Cref{d:rds}(ii) holds.  Thus, it remains to show 
	\begin{equation}\label{e:srtrees}
		V(\bs\pi_{G_{1,1}})\otimes V(\bs\pi_G) \quad\text{is simple}. 
	\end{equation}
	In fact, \eqref{e:kkopineq} implies
	\begin{equation*}
		\lie d(G_{1,1},G)\leq \lie d(G_{1,1},G_1) + \lie d(G_{1,1},G'').
	\end{equation*}
	The second summand is zero by \eqref{e:srtreenov1} and \Cref{kkop}(a), while the first is zero since $G_{1,1}$ satisfies \Cref{d:rds}(iii) when regarded as a subgraph of $G_1$. Thus, $\lie d(G_{1,1},G)=0$ and we are done after one final application of \Cref{kkop}(a). 
\end{proof}

\begin{rem}
	As in \Cref{r:treesr}, the only reason we needed to assume $\lie g$ is of type $A$ in \Cref{t:srtree} is because we used \Cref{l:3aline} (via \Cref{c:eas2}) for proving \eqref{e:srtrees} for the base of the inductive argument. Thus, the theorem actually holds for $\lie g$ of any type as long as the conclusion of \Cref{l:3aline} holds for all the pairs of linking vertices of the given $\mathcal G$-tree. \endd
\end{rem}

Let us finish this section with an example of a strongly real module whose $q$-factorization graph is one step away from being in the realm of \Cref{t:srtree}.

\begin{ex}\label{ex:srnst}
	Suppose that $\lie g$ is of type $A_n, n\ge 3$, and let $\bs\pi =\bs\omega_{1,2,3}\, \bs\omega_{3,6,3}\, \bs\omega_{2,9,3}$. One easily checks $V(\bs\pi)$ is not a snake module, while $G(\bs\pi)$ is 
	
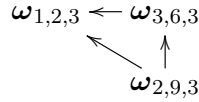
\begin{equation*}
		\xymatrix@R-1pc@C-1pc{ \bs\omega_{1,2,3} & \ar[l]\bs\omega_{3,6,3}\\ & \ar[u]\ar[ul]\bs\omega_{2,9,3}}
	\end{equation*}
	Notice every two-vertex subgraph is a snake tree, but $G(\bs\pi)$ cannot be realized as a snake tree, so \Cref{t:srtree} cannot be used. In fact, one can also check that no pseudo $q$-factorization graph over $\bs\pi$ can be realized as a snake tree.
	
    We claim that $H=G(\bs\omega_{3,6,3})$ is an rds for $G(\bs\pi)$. Since $H^c$ is a tree, it follows that $V(\bs\pi)$ is strongly real. To prove the claim first note, we have already noted that $H$ determines a real cut.  \Cref{c:hlwpqf} implies 
	\begin{equation}\label{e:srnsthw}
		V(\bs\omega_{2,11})\otimes V(\bs\omega_{2,8,2}\,\bs\omega_{1,2,3})\quad\textrm{and}\quad V(\bs\omega_{2,11})\otimes V(\bs\omega_{3,6,3}) \quad\text{are highest-$\ell$-weight}.
	\end{equation}
	Also,  as we will check below, it follows from \Cref{t:3lineprime} that 
	\begin{equation}\label{e:srnsts}
		V(\bs\omega_{2,8,2}\,\bs\omega_{1,2,3})\otimes V(\bs\omega_{3,6,3}) \quad\text{is simple.}
	\end{equation} 
	The checking that condition (ii) of \Cref{d:rds} is satisfied now follows since \eqref{e:srnsthw} and \eqref{e:srnsts}, together with \Cref{c:hlwquot}, imply
	\begin{equation*}
		V(\bs\pi_{H^c})\otimes V(\bs\pi_H) = V(\bs\omega_{2,9,3}\,\bs\omega_{1,2,3})\otimes V(\bs\omega_{3,6,3})\quad \text{is highest-$\ell$-weight.}
	\end{equation*}
	
	The above argument also implies   
	\begin{equation*}
		V(\bs\pi_{H^c})\otimes V(\bs\pi_H)\otimes  V(\bs\pi_H) \quad \text{is highest-$\ell$-weight}
	\end{equation*}
	and, hence, so is $V(\bs\pi)\otimes V(\bs\pi_H)$. Therefore, in order to check that condition (iii) of \Cref{d:rds} holds, it suffices to prove that 
	\begin{equation}\label{e:srnsthw'}
		V(\bs\pi_H)\otimes V(\bs\pi) \quad \text{is highest-$\ell$-weight.}
	\end{equation}
	For doing that, note that the modules $V(\bs\omega_{3,6,3})\,\otimes V(\bs\omega_{1,2,3})$ and $V(\bs\omega_{2,9,3}\,\bs\omega_{3,6,3})\otimes V(\bs\omega_{1,2,3})$ are both highest-$\ell$-weight by \Cref{c:hlwquot}. This, together with \eqref{e:srnsts} and \Cref{c:hlwquot} implies
	\begin{equation*}
		V(\bs\omega_{3,6,3})\otimes V(\bs\omega_{2,9,3}\,\bs\omega_{3,6,3})\otimes V(\bs\omega_{1,2,3})
	\end{equation*}
	is highest-$\ell$-weight and projects onto $V(\bs\pi_H)\otimes V(\bs\pi)$, thus proving \eqref{e:srnsthw'}.  We have checked that neither the sink nor the source of this graph are rds. Thus, this is an example of a strongly real graph which do not admit an rds-quochain whose rds are extremal vertices. 
	 Let us check the assumptions of \Cref{t:3lineprime} are satisfied so that \eqref{e:srnsts} indeed follows from it. In the notation of the theorem, let $i_j=3$, so $i_{j'}=2, m_j = 4, m_{j'}=6, I_j = [1,3], I_{j'}=[1,2]$. Then, $i_{j'}\in I_{j'}\subseteq I_j$, which also implies $m_{j'}\in \mathscr{R}_{i,i_{j'},I_j}^{r,r_{j'}}$. The third condition reads $6-4+4\in\mathscr R_{1,2,[1,3]}^{3,2} = \{4,6\}$, while \eqref{2exracondm} reads $4+3\le 6+2+1$. Finally, we need to check the graph is an alternating line. For doing this, it remains to check that there is no arrow between the vertices corresponding to $\bs\omega_{2,8,2}$ and $\bs\omega_{3,6,3}$,  i.e., we need to show $2 = m_{j'}-m_j\notin\mathscr R_{2,3}^{2,3}$. But $\mathscr R_{2,3}^{2,3}=\{4,6,8\}$ by \eqref{e:typeAR}.\endd	
\end{ex}
			
\subsection{RDS-Quochains via KKOP Invariants}\label{ss:kkkords}
We now discuss a method for constructing real modules afforded by $\mathcal G$-trees for which the set $\mathcal A_\mathcal G$ are related to KKOP invariants whose values are at most $1$. This assumption allows us to use \Cref{c:kkop} in order to ensure that  \Cref{t:isreal}(iii) holds, or equivalently, that \Cref{d:rds}(iii) is satisfied. 

\begin{thm}\label{t:kkoprds}
	Let $\mathcal G = G_1,\dots,G_l$ be a sequence of pseudo $q$-factorization graphs such that $G_k$ is real for all $1\le k\le l$. If $G=G_1\otimes \cdots\otimes G_l$ is a $\mathcal G$-tree  and 
	\begin{equation*}
		\lie d (V(\bs\pi_{G_k}),V(\bs\pi_{G_m}))\leq 1 \quad\text{for all}\quad 1\leq m,k\leq l,
	\end{equation*}
	then  $V(\bs\pi_G)$ is real. More precisely, if the sequence is numbered so that ${\rm val}_{\bar G_{k-1}}(G_k)=1$ for all $1\le k<l$, then $\mathcal G$ is a weak rds-quochain for $G$. 
\end{thm}

\begin{proof}
	We proceed by induction on $l$ which clearly begins for $l=1$. Thus, assume $l>1$. 	
	Since $\mathcal G$ is a $\mathcal G$-tree, there exists $1\leq k\leq l$ such ${\rm val} (G_k)=1$. Without loss of generality, assume $k=1$ and that $G_1$ is linked to $G_2$. 
	
	By the induction hypothesis $G'=G_2\otimes\cdots \otimes G_k$ is real and \Cref{kkop}(a) implies 
	\begin{equation*}
		\lie d (V(\bs\pi_{G_1}),V(\bs\pi_{G_m})) = 0\ \ \text{if}\ \ m\ge 3.
	\end{equation*}
	It then follows from \eqref{e:kkopineq} that 
	\begin{equation*}
		\lie d (V(\bs\pi_{G_1}),V(\bs\pi_{G'}))\leq \lie d (V(\bs\pi_{G_1}),V(\bs\pi_{G_2}))\leq 1.
	\end{equation*}
	An application of \Cref{l:rdsbysmallkkop} shows $G_1$ is an rds and, hence, the second statement follows from the induction assumption.
\end{proof}

Let us give an example where we can use a combination of Theorems \ref{t:srtree} and \ref{t:kkoprds} to obtain a strong rds-quochain, but neither of these theorems can produce such quochain directly.

\begin{ex}\label{ex:srnstkkop}
	Let $\bs\pi$ be as in \Cref{ex:srnst} and, given $m\in\mathbb Z$, denote by $\tau_m$ the group automorphism of $\mathcal P$ determined by $\bs\omega_{i,a}\mapsto \bs\omega_{i,a+m}$ for all $i\in I,a\in\mathbb Z$. Let $\bs\pi_k = \tau_{14(k-1)}\bs\pi$ for $1\le k\le l$ and consider $\bs\varpi = \bs\pi_1\bs\pi_2\cdots\bs\pi_l$.
	One easily checks $G(\bs\varpi)$ is 
	\begin{equation*}
		\xymatrix@R-1pc@C-1pc{ 
			\bs\omega_{1,2,3} & \ar[l]\bs\omega_{3,6,3} & \bs\omega_{1,16,3}\ar[dl] & \bs\omega_{3,20,3}\ar[l]&\ar[dl]\cdots & \bs\omega_{1,14(l-1)+2,3}\ar[dl] & \bs\omega_{3,14(l-1)+6,3}\ar[l]\\ 
							  & \ar[u]\ar[ul]\bs\omega_{2,9,3} & & \ar[u]\ar[ul]\bs\omega_{2,23,3}&\cdots &&\ar[u]\ar[ul]\bs\omega_{2,14(l-1)+9,3}}
	\end{equation*}
	Setting $G_k=G(\bs\pi_k)$ and $\mathcal G = G_1,\dots, G_l$, it follows from \Cref{ex:srnst} that $G_k$ is (strongly) real and, moreover,  since there are no arrows linking $G_k$ and $G_{k'}$ if $|k-k'|>1$, 
	\begin{equation*}
		|k-k'|>1 \quad\Rightarrow\quad \lie d(V(\bs\pi_k), V(\bs\pi_{k'}))=0.
	\end{equation*}
	Furthermore, \Cref{kkop} implies
	\begin{equation}\label{e:srnstkkop}
		\lie d(V(\bs\pi_k), V(\bs\pi_{k+1}))\leq \lie d(V(\bs\pi_k), V(\bs\omega_{1,14k+2,3}))\leq  \lie d(V(\bs\omega_{2,14(k-1)+9,3}), V(\bs\omega_{1,14k+2,3})).
	\end{equation}
	Finally, one easily checks $V(\bs\omega_{2,14(k-1)+9,3}\,\bs\omega_{1,14k+2,3})$ is a prime snake module and a standard iterated application of \Cref{l:naoi} implies that the above KKOP invariant is at most $1$. Thus, \Cref{t:kkoprds} implies $\mathcal G$ is a weak rds-quochain for $G(\bs\varpi)$. 
	
	Note that \Cref{t:srtree} cannot be used to describe rds quochains for $G$. Indeed, the realization of $G$ as a $\mathcal G$-tree is not well-based since the tails and heads of the arrows in $\mathcal A_{\mathcal G}$ are not rds of the corresponding parts of the multicut by \Cref{ex:srnst}. Still, using the above described weak rds quochain and a mix of the ideas of the proofs of Theorems \ref{t:srtree} and \ref{t:kkoprds}, it can be shown that $G$ is strongly real as follows.  
	
	\Cref{ex:srnst} implies the singleton-sequence 
	\begin{equation*}
		\mathcal G_k = \bs\omega_{3,14(k-1)+6,3},\ \bs\omega_{1,14(k-1)+2,3},\ \bs\omega_{2,14(k-1)+9,3}
	\end{equation*}
	is an rds-quochain for $G_k$. Let us check the concatenation $\mathcal G_1*\cdots*\mathcal G_l$ is an rds-quochain for $G$. As usual, we proceed by induction on $l$, which starts when $l=1$ by \Cref{ex:srnst}.
	Let $l>1$, $\mathcal G' = G_2,\dots,G_l$, and $G'=G_2\otimes\cdots\otimes G_l$. By the induction assumption, $\mathcal G_2*\cdots*\mathcal G_l$ is an rds-quochain for $G'$. Let also 
	\begin{equation*}
		G_{k,1} = G(\bs\omega_{3,14(k-1)+6,3}), \quad G_{k,2} = G(\bs\omega_{1,14(k-1)+2,3}), \quad G_{k,3} = G(\bs\omega_{2,14(k-1)+9,3}),
	\end{equation*} 
	so that $\mathcal G_k = G_{k,1}*G_{k,2}*G_{k,3}$. 	We begin by checking that 
	\begin{equation}\label{e:srnstkkop3}
		G_{1,3} \ \ \text{is an rds for}\ \ G_{1,3}\otimes G'.
	\end{equation}
	But $G_{1,3}\otimes G'$ is an $\mathcal H$-tree where $\mathcal H=H_1,H_2$ with $H_1 = G_{1,3}$ and $H_2=G'$. We already know $G_{1,3}$ determines a real cut and, since it is a bottom subgraph,
	\begin{equation*}
		V(\bs\pi_{G'})\otimes V(\bs\omega_{2,14(k-1)+9,3}) \quad\text{is highest-$\ell$-weight}.
	\end{equation*}
	Moreover, \eqref{e:srnstkkop} can be used as before to conclude $\lie d(G_{1,3},G'\otimes G_{1,3})\le 1$, so \Cref{t:kkoprds} concludes the checking of \eqref{e:srnstkkop3}. 
	
	Note that, if we set $H = G_{1,1}\otimes G_{1,2}$, then $H^c = G_{1,3}\otimes G'$ and the cut $(H,H^c)$ is real and formed by extremal subgraphs. However, it does not provide a realization of $G$ an $\mathcal H$-tree with $H_1=H$ and $H_2=H^c$ since $\mathcal A_{\mathcal H}$ has two arrows. Thus, we cannot use our results about generalized trees directly. However, if we set $H_1= G_{1,2}$ and $H_2= G_{1,3}\otimes G'$, then $\mathcal H= H_1,H_2$ is an $\mathcal H$-tree realization of $G\setminus G_{1,1}$ since $\mathcal A_{\mathcal H}$ has a single arrow, which links $G_{1,2}$ to $G_{1,3}$. Since $G_{1,2}$ is and rds-base for itself and \eqref{e:srnstkkop3} implies $G_{1,3}$ is an rds-base for $H_2$, it follows that we can use \Cref{t:srtree} to conclude 
	\begin{equation}\label{e:srnstkkop2}
		G_{1,2} \ \ \text{is an rds for}\ \ G\setminus G_{1,1}.
	\end{equation}
	It remains to check
	\begin{equation}\label{e:srnstkkop1}
		G_{1,1} \ \ \text{is an rds for}\ \ G.
	\end{equation}
	No multicut having $G_{1,1}$ as a part gives rise to a realization of $G$ as a generalized tree. So, we cannot use our theorems, but the ideas of the proof of \Cref{t:srtree} can be used to check the above. Namely, since we have just shown $G\setminus G_{1,1}$ is real, it determines a real cut. On the other hand, $G_{1,1}$ is an rds for $G_1$ and, hence, $\lie d(G_{1,1},G_1)=0$. This proves \Cref{d:rds}(ii) holds and implies 
	\begin{equation*}
		\lie d(G_{1,1},G) \le  \lie d(G_{1,1},G') + \lie d(G_{1,1},G_1) = \lie d(G_{1,1},G'). 
	\end{equation*}
	But the latter is zero since there are no arrows linking $G_{1,1}$ and $G'$. This completes the checking that $\mathcal G_1*\cdots*\mathcal G_l$ is an rds-quochain formed by singletons and, hence, $G$ is strongly real.	\endd
\end{ex}

Let us end this section with an example of an application of \Cref{t:kkoprds} with a $\mathcal G$-tree whose parts are not strongly real. 

\begin{ex}\label{ex:realbutnotstrongkkop}
	Let $\bs\pi_1 =2_01_33_32_6^2$ as in \Cref{ex:realbutnotstrong}, in type $A_3$,  and set $\bs\pi_2 =2_{-4}1_{-7}3_{-7}2_{-10}^2$. Note that, up to a uniform shift of the centers of the $q$-strings,  $\bs\pi_2 = \bs\pi_1^{\kappa}$.  In particular,  $G= G(\bs\pi_1\bs\pi_2)= G_f(\bs\pi_1\bs\pi_2)$ is given by 
	
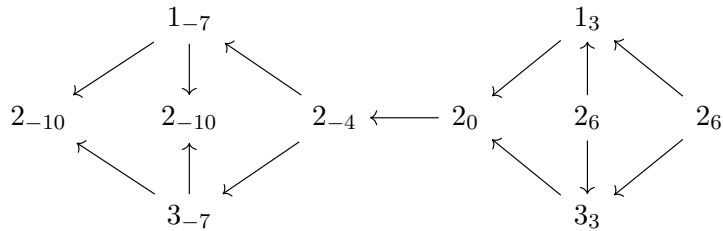
\begin{equation*}
		\begin{tikzcd}
			& 1_{-7}\ar[dl]\ar[d] & & & 1_3\ar[dl]& \\ 
			2_{-10}&2_{-10}& 2_{-4}\ar[ul]\ar[dl]& 2_0\ar[l] & 2_6\ar[u]\ar[d]&2_6\ar[ul]\ar[dl]\\ 
			& 3_{-7}\ar[ul]\ar[u] & & & 3_3\ar[ul]&  
		\end{tikzcd}
	\end{equation*}
	which is realizable as a $\mathcal G$-tree with $\mathcal G = G_1, G_2$, and $G_k=G(\bs\pi_k), k\in\{1,2\}$. Moreover, the usual multiple application of \Cref{kkop}, alongside  \Cref{l:naoi}, implies 
	\begin{equation*}
		\lie d (V(\bs\pi_1),V(\bs\pi_2))\leq \lie d (V(2_{-4}),V(2_0))= 1.
	\end{equation*}
	This, together with the fact that $G_k$ is real by the analysis of \Cref{ex:realbutnotstrong}, shows the assumptions of \Cref{t:kkoprds} are satisfied and, hence, $\mathcal G$ is a weak rds-quochain for $G$. However, as seen in \Cref{ex:realbutnotstrong}, the rds $G_1$ and $G_2$ are not strongly real. \endd
\end{ex}

\subsection{Comments on Primality and Cluster Algebras}\label{ss:cluster}
Note the graphs in Examples \ref{ex:prst}, \ref{ex:srnst}, and \ref{ex:srnstkkop} are the actual $q$-factorization graphs and they are totally ordered. Hence, the corresponding modules are also prime by \cite[Theorem 3.5.5]{ms:to}. 
Let us check that the one in \Cref{ex:nextglue} is also prime. 

\begin{ex}\label{ex:nextgluep}
	Let $\bs\pi$ and $G$ be as in \Cref{ex:nextglue} and let us check $V(\bs\pi)$ is prime. Suppose $\bs\pi_1,\bs\pi_2\in\mathcal P^+$ are such that 
	\begin{equation*}
		V(\bs\pi)\cong V(\bs\pi_1)\otimes V(\bs\pi_2).
	\end{equation*} 
	As seen in \cite{ms:to}, it suffices to consider the case that $(\bs\pi_1,\bs\pi_2)=(\bs\pi_H,\bs\pi_{H^c})$ with $H$ and $H^c$ being connected subgraphs of $G$.  
	Without loss of generality assume $\bs\omega_{1,7}|\bs\pi_1$. If $\bs\omega_{2,4}|\bs\pi_2$, letting $J=[1,2]$ and using \eqref{e:defredset}, it immediately follows that $V(\bs\pi_1)_J\otimes V(\bs\pi_2)_J$ is reducible, contradicting \eqref{e:sJs}. Therefore $\bs\omega_{1,7}\bs\omega_{2,4}|\bs\pi_1$. A similar argument with $J=[2,3]$ imply $\bs\omega_{3,9,3}|\bs\pi_1$ and, then, a similar argument with $J=[3,4]$ shows $\bs\omega_{4,0}|\bs\pi_1$. Thus, $\bs\pi_2=\bs 1$, which completes the checking.\endd
\end{ex}

Prime real modules are related to cluster algebras and it would be interesting to study how the examples above  appear in that context. Apart from specific classes of modules, such as snake modules \cite{DLL:snakes}, HL-modules \cite{HL:cluster, bc:tphl,DS:hlade},  and their generalizations \cite{DGL:hlgen}, as far as we know, there is no explicit description of Drinfeld polynomials whose associated irreducible representation is prime and real and  correspond to cluster variables. We leave this discussion for a future publication. 

In \cite{ms:3tree}, a complete description of the prime modules for type $A$ whose $q$-factorization graphs have three vertices was given. The most difficult case is that of graphs which are alternating lines  (\Cref{t:3lineprime}). It would also be interesting to consider the case of snake trees formed by three snakes, as well as a triangle of snakes. In \cite{BC:nil}, a particular class of snake trees is studied and the authors characterize the prime ones.

\subsection{Are There Real Modules With No RDS?}\label{ss:nords}
Frustratingly and intriguingly enough, we have not found an example of pseudo $q$-factorization associated to a real module which admits no rds (cf. \cite[Remark 7.7]{LM:sqirred}). Let us end this paper with a discussion related to this from the perspective of \Cref{d:weakrdsq}. 

Namely, this definition can be used to create a stratification of the class of real modules by a level of ``complexity'' measured by the length of weak rds-quochains. For simplicity, let us work with actual $q$-factorization graphs only. Given $\bs\pi\in\mathcal P^+$, define the rds-quochain index of $\bs\pi$ by
\begin{gather*}
	Q(\bs\pi) = \max\{l: \exists\ \text{a weak rds quochain of length $l$ for $G(\bs\pi)$}\},
\end{gather*}
where we understand $Q(\bs\pi)=0$ if there does not exist a weak rds quochain for $G(\bs\pi)$. Hence,
\begin{equation*}
	Q(\bs\pi) = 0 \quad\Leftrightarrow\quad V(\bs\pi) \ \ \text{is imaginary}.
\end{equation*}
Then, the ``reality index'' of $\bs\pi$ can be defined as
\begin{gather*}
	R(\bs\pi) = \#\mathcal V_{G(\bs\pi)} - Q(\bs\pi),
\end{gather*}
so 
\begin{equation*}
	R(\bs\pi) = \#\mathcal V_{G(\bs\pi)} \quad\Leftrightarrow\quad V(\bs\pi) \ \ \text{is imaginary}
\end{equation*}
and
\begin{equation*}
	R(\bs\pi)=0 \quad\Leftrightarrow\quad G(\bs\pi)\ \ \text{is strongly real}.
\end{equation*}
The real modules which are not strongly real  satisfy 
\begin{equation*}
	0<R(\bs\pi)<\#\mathcal V_{G(\bs\pi)}.
\end{equation*} This number measures how far a module is from being strongly real. Strongly real modules should be regarded as the ``simplest'' real modules as their Drinfeld polynomials can be built by a sequence of KR factors such that all the intermediate modules are real.

The modules with $Q(\bs\pi)=1$ should be regarded as  ``basic'' real modules, since they cannot be built from two smaller real modules. KR modules are basic and strongly real. They are the only basic real modules which are strongly real. 
The question in the title of this subsection is then equivalent to: are there basic real modules beyond KR modules? 
Such modules would have
\begin{equation*}
	R(\bs\pi) =  \#\mathcal V_{G(\bs\pi)} -1 >0. 
\end{equation*}
If KR modules are the only basic real modules, then
\begin{equation*}
	Q(\mathcal P^+\setminus\mathcal{KR}) \subseteq \mathbb Z_{\ge 0}\setminus\{1\}
\end{equation*}
and all real modules with at least two $q$-factors can be constructed as a subquotient of two real modules whose Drinfeld polynomials have strictly smaller degree. This is rather intriguing and there should be a deeper explanation if that is indeed the case.

Examples of real modules which are not strongly real are given in Examples \ref{ex:realbutnotstrong} and \ref{ex:realbutnotstrongkkop}. However, they admit proper rds and, hence, are not basic. In fact, in the case of \Cref{ex:realbutnotstrong}, we have shown $Q(\bs\pi) = 4 = \#\mathcal V_{G(\bs\pi)} -1$, so $R(\bs\pi)=1$ and the module is in the stratum determined by $R(\bs\pi)$ immediately above the stratum of strongly real modules.  This, together with \Cref{t:kkoprds}, implies that, for $\bs\pi$ as in \Cref{ex:realbutnotstrongkkop},  we have $R(\bs\pi)\le 2$. A careful analysis such as the one we have performed in \Cref{ex:realbutnotstrong} should show that equality holds. We shall leave a more systematic study of such considerations for the future.


\begin{thebibliography}{10}
	
	\bibitem{bc:tphl}
	M. Brito and V. Chari, {\em Tensor products and $q$-characters of HL-modules and monoidal categorifications}, \href{http://www.numdam.org/articles/10.5802/jep.101}{Journal de l’École polytechnique — Mathématiques, 6 (2019), 581-619}
	
		
	\bibitem{BC:nil}
	M. Brito and V. Chari, {\em Alternating snake modules and a determinantal formula}, \href{https://doi.org/10.48550/arXiv.2412.03750}{arXiv:2412.03750}.
	
	\bibitem{bms:tof}
	M. Brito, A. Moura, and C. Silva, {\em Totally ordered pseudo $q$-factorization graphs and prime factorization}, \href{https://arxiv.org/abs/2410.01519}{to appear in Arkiv f\"or Matematik}.
	
	\bibitem{cdf:grass}
	W. Chang, B. Duan, C. Fraser, J-R. Li, {\em Quantum affine algebras and Grassmannians}, \href{https://doi.org/10.1007/s00209-020-02496-7}{Math. Z. 296, 1539–1583 (2020)}. 
	
	\bibitem{DGL:hlgen}
	B. Duan, J. Guo, Y. Luo, {\em Generalized Hernandez-Leclerc modules and cluster algebras}, \href{https://doi.org/10.1142/S0219498826501094}{Journal of Algebra and Its Applications} .
	
	\bibitem{DLL:snakes} 
	B. Duan, J.-R. Li and Y.-F. Luo, {\em Cluster algebras and snake modules}, \href{https://doi.org/10.1016/j.jalgebra.2018.10.027}{Journal of Algebra, Volume 519, 2019, 325-377}.
	
	\bibitem{DS:hlade}
	B. Duan and R. Schiffler, {\em Real simple modules over simply-laced quantum affine algebras and categorifications of cluster algebras}, \href{https://arxiv.org/abs/2305.08715}{arXiv:2305.08715
	}
	
	\bibitem{eali:trop}
	N. Early and J-R. Li, {\em Tropical geometry, quantum affine algebras, and scattering amplitudes}, \href{https://iopscience.iop.org/article/10.1088/1751-8121/ad909b}{J. Phys. A: Math. Theor. 57 (2024), 495201}.
	
	
	\bibitem{gumi}
	M. Gurevich and A. Mínguez, {\em Cyclic representations of general linear p-adic groups}, \href{https://doi.org/10.1016/j.jalgebra.2021.05.013}{J. Alg. 585 (2021), 25--35}.
	
	\bibitem{HL:cluster}
	D. Hernandez and B. Leclerc, {\em Cluster algebras and quantum affine algebras}, \href{https://doi.org/10.1215/00127094-2010-040}{Duke Math. J. 154 (2) 265 - 341, 15 August 2010}.
	
	\bibitem{hl:KR}
	D. Hernandez and B. Leclerc, {\em A cluster algebra approach to $q$-characters of Kirillov-Reshetikhin modules} \href{https://doi.org/10.4171/jems/609}{J. Eur. Math. Soc. 18 (2016), no. 5, 1113–1159}
	

	
	\bibitem{kkop:mc1}
	M. Kashiwara, M. Kim, S. Oh and E. Park, {\em Monoidal categorification and quantum affine algebras}, \href{https://doi.org/10.1112/S0010437X20007137}{Compositio Mathematica. 156(5):1039-1077, 2020}. 
	
	\bibitem{KKOP:cluster} 
	M. Kashiwara, M. Kim, S.-J. Oh and E. Park, {\em Cluster algebra structures on module categories over quantum affine algebras}, \href{https://doi.org/10.1112/plms.12428}{Proc. London Math. Soc., (2022) 124: 301-372}. 
	
	\bibitem{kkop:pbw}
	M. Kashiwara, M. Kim, S-J. Oh, E. Park, {\em PBW theory for quantum affine algebras}. \href{https://doi.org/10.4171/JEMS/1323}{J. Eur. Math. Soc. (2023), published online first}.
	
	\bibitem{kkop:mc2}
	M. Kashiwara, M. Kim, S-J. Oh, E. Park, {\em Monoidal categorification and quantum affine algebras II}.\href{https://doi.org/10.1007/s00222-024-01249-1}{Invent. math. 236, 837–924 (2024)}. 
	
	
	\bibitem{LM:sqirred} 
	E. Lapid and A. M\'inguez, {\em Geometric conditions for $\square$-irreducibility of certain representations of the general linear group over a non-archimedean local field}, \href{https://doi.org/10.1016/j.aim.2018.09.027}{Adv. Math. 339 (2018), 113–190}.
	
	
	\bibitem{lec:im}
	B. Leclerc, {\em Imaginary vectors in the dual canonical basis of $U_q(n)$}, Transformation Groups, 8:95–104, 2002.
	
	\bibitem{Moura}
	A.~Moura, {\em An introduction to finite-dimensional representations of classical and quantum affine algebras}, lecture notes published in \href{https://www.famaf.unc.edu.ar/documents/887/BMat59.pdf}{Trabajos de matem\'atica s\'erie B 59}, Publicaciones de la FaMAF - Universidad Nacional de C\'ordoba, 2011. 
	
	\bibitem{ms:to} 
	A. Moura and C. Silva, {\em On the primality of totally ordered $q$-factorization graphs}, \href{http://dx.doi.org/10.4153/S0008414X23000160}{Canad. J. Math., 76  (2024), 594–637}.
	
	\bibitem{ms:3tree} 
	A. Moura and C. Silva, {\em Three-vertex prime graphs and reality of trees}, \href{https://doi.org/10.1080/00927872.2023.2196345}{Comm. in Algebra, 51 (2023), 4054–4090}.
	
	\bibitem{my:pathB}
	E. Mukhin and C. Young, {\em Path description of type B $q$-characters}, \href{https://doi.org/10.1007/s00029-011-0083-x}{Adv. Math. 231:1119–1150, 2012}.
	
	\bibitem{muyou:tsystem}
	E. Mukhin and C. Young, {\em Extended T-systems}, \href{https://doi.org/10.1007/s00029-011-0083-x}{Selecta Mathematica {\bf 18} (2012), 591--631}.
	
	\bibitem{naoi:Tsys}
	K. Naoi, {\em Strong duality data of Type A and extended T-Systems}, \href{https://doi.org/10.1007/s00031-024-09860-5}{Transformation Groups (2024)}. 
	
	\bibitem{ohscr:simptens}
	S. J. Oh and T. Scrimschaw, {\em Simplicity of tensor products of Kirillov-Reshetikhin modules: nonexceptional affine and $G$ types}, \href{https://arxiv.org/abs/1910.10347}{arXiv:1910.10347}.
	
	\bibitem{Qin}
	F. Qin, {\em Triangular bases in quantum cluster algebras and monoidal categorification conjectures}, Duke Math. J., 166 (12) (2017), 2337-2442.
	
\end{thebibliography}
\end{document}